\documentclass[number,citesort,seceqn,dvips]{arxbj}
\usepackage{graphicx}

\aid{0}
\volume{17}
\issue{4}
\pubyear{2011}
\firstpage{1285}
\lastpage{1326}
\doi{10.3150/10-BEJ321}

\makeatletter
\newcommand{\Nset}{\mathbb{N}}
\newcommand{\Rset}{\mathbb{R}}
\renewcommand{\epsilon}{\varepsilon}
\newcommand{\G}{\mathcal{G}}

\newtheorem{theo}{Theorem}[section]
\newtheorem{prop}{Proposition}[section]
\newremark{rem}{Remark}[section]
\newtheorem{lem}{Lemma}[section]

\makeatother

\begin{document}
\begin{frontmatter}

\title{Probabilistic sampling of finite renewal~processes}
\runtitle{Probabilistic sampling of finite renewal~processes}

\begin{aug}
\author[1]{\fnms{Nelson} \snm{Antunes}\thanksref{1,t1}\ead[label=e1]{nantunes@ualg.pt}}
\and
\author[2]{\fnms{Vladas} \snm{Pipiras}\thanksref{2,t2}\corref{}\ead[label=e2]{pipiras@email.unc.edu}}
\runauthor{N. Antunes and V. Pipiras}
\address[1]{University of Algarve/CEMAT, Faro, Portugal and INRIA Paris -- Rocquencourt, Paris, France. \printead{e1}}
\address[2]{Instituto Superior T\'{e}cnico/CEMAT, Lisboa, Portugal and University of North Carolina, Chapel Hill, NC, USA. \printead{e2}}
\thankstext[*]{t1}{Current address: Campus de Gambelas, 8900 Faro, Portugal.}
\thankstext[**]{t2}{Current address: Av. Rovisco Pais, 1049-001 Lisboa, Portugal.}
\end{aug}

\received{\smonth{1} \syear{2010}}
\revised{\smonth{8} \syear{2010}}

\begin{abstract}
Consider a finite renewal process in the sense that interrenewal times
are positive i.i.d. variables and the total number of renewals is a
random variable, independent of interrenewal times. A finite point
process can be obtained by probabilistic sampling of the finite renewal
process, where each renewal is sampled with a fixed probability and
independently of other renewals. The problem addressed in this work
concerns statistical inference of the original distributions of the
total number of renewals and interrenewal times from a sample of i.i.d.
finite point processes obtained by sampling finite renewal processes.
This problem is motivated by  traffic measurements in the Internet in
order to characterize flows of packets (which can be seen as finite
renewal processes) and where  the use of packet sampling is becoming
prevalent due to increasing link speeds and limited storage and
processing capacities.
\end{abstract}


\end{frontmatter}

\section{Introduction}
\label{s:intro}

\subsection{Motivation}

The statistical and probabilistic problems considered in this work are
motivated by questions arising from data traffic analysis in modern
communication networks such as the Internet. Over these networks,
information is sent in the form of packets, and packets are grouped
into flows. (A flow corresponds to a group of packets sharing common
characteristics. Ideally, a flow can be thought of as a set of packets
that arises in the network through a remote terminal session or a web
page download.) Each packet carries information about the flow it
belongs to and also whether it is the first or the last packet in the
flow. Examining each packet then allows all flows to be reconstructed.
Knowing the structure of flows (flow level characteristics) helps
network operators and networking researchers to understand and
discover characteristics of data traffic.

A key difficulty with capturing each packet is that this rapidly leads
 to a huge amount of data to store and  analyze. For
instance, capturing all traffic for a few hours on a gigabit/sec. link
at a~medium load level yields several hundred  gigabytes of data. A way
to reduce the volume of data is by sampling packets. One of the
simplest approaches is probabilistic sampling, where each packet,
independently of the others, is captured and analyzed with a fixed
probability $q$. The basic problem, known as the \textit{inversion
problem} in the network traffic literature, is then to  deduce the
structure of the original flows from sampled packets.

This problem has recently attracted much attention  in the networking
community, with the focus almost exclusively on inference of the
original distribution of flow sizes (total number of packets);  see
\cite{duffieldlundthorip2006,hohnveitch2006,yangmichailidis2007}.
In this work, we take a closer look at statistical properties of the
previously considered estimator of the original distribution of flow
sizes and are also interested in inference of interarrival times
between packets from sampled data. This distribution allows  traffic
burstiness to be characterized and leads to other major characteristics
such as the duration of a flow.

\subsection{Statement of problem and goals}

\begin{figure}[b]

\includegraphics{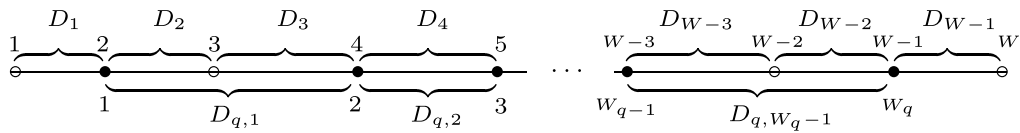}

\caption{Finite renewal process and sampled finite point process.}
\label{FRP-SSP}
\end{figure}

From a mathematical standpoint, flows and their probabilistic sampling
can be described as follows. As suggested by data traffic measurements
(e.g., \cite{hohnveitchabry2003}), a flow can be modeled as a finite
renewal process, where the total number of renewals $W$ is random and
the interrenewal times $D_i$, $i= 1,\ldots,W-1$, are positive i.i.d.
random variables, independent of $W$. Suppose that the finite renewal
process is sampled probabilistically with a fixed probability $q$ to
obtain a sampled finite point process. (As shown below, the resulting
sampled point process is generally not a finite renewal process.) The
sampled number of renewals (which could be zero) and sampled
interrenewal times (which are available only when the number of sampled
renewals is greater than $1$) will be denoted by $W_q$ and $D_{q,i}$,
$i=1,\ldots,W_q-1$, respectively. We illustrate this notation in Figure
\ref{FRP-SSP}, where all circles (both empty and full) correspond to
the finite renewal process, and full circles correspond to the sampled
finite point process.

Given $N$ i.i.d. copies of a sampled finite point process (some of
which are empty when the number of sampled renewals is zero), one of
our main goals is to infer the original distribution $F_D$ of
interrenewal times $D_i$. The focus is on nonparametric inference of
$F_D$, statistical properties of the resulting estimator and its
performance in simulations. The estimator of $F_D$ involves a~nonparametric estimator of the original distribution $f_W$ of the total
number of renewals $W$, which has  previously been considered in the
network traffic context \cite{hohnveitch2006}.
Statistical properties of the estimator of $f_W$ will also be studied
here for the first time.

\subsection{Inference procedure}

To be able to make inference about $F_D$, we first need to relate it to
$F_{D_{q,i}|s}$, where the latter indicates the conditional
distribution function of $D_{q,i}$ given $W_q=s$ (with $i<s$). As part
of this work (see Theorem \ref{t:LT_sit} in Section
\ref{s:renewal-to-sampled}), we show that
%
\begin{equation}\label{e:Dq-D-Wknown}
    F_{D_{q,i}|s} = \sum_{m=1}^\infty A_{s,m} F_D^{*m},
\end{equation}
where $*m$ denotes the $m$th convolution. In particular, the right-hand
side of (\ref{e:Dq-D-Wknown}), and hence $F_{D_{q,i}|s}$, does not
depend $i$. We also note for later reference that the definition of the
sequence $A_s = \{A_{s,m}\}_{m\in \Nset}$ involves the distribution
$f_W$. Somewhat independently of the main objectives above, we examine
(\ref{e:Dq-D-Wknown}) for several underlying distributions $f_W$, such
as geometric and heavy-tailed, that arise in network traffic studies
(e.g., \cite{chabchoubfrickerguilleminrobert2010}). For example, in the
geometric case, the distribution $F_{D_{q,i}|s}$ does not depend on
$s$, although this seems to be an exception to the general
rule. 
We also provide a result similar to (\ref{e:Dq-D-Wknown}) for a
multidimensional vector $(D_{q,1},\ldots,D_{q,n})$. This allows
(conditional) dependencies to be examined among sampled interrenewal
times. For example, when the distribution $f_W$ is geometric, the
sampled interrenewal times turn out to be independent, in which case
the sampled finite point process is a finite renewal process. We should
also note that forms of conditioning other than on $W_q=s$ are
possible, such as $W_q\geq s$; these will be discussed briefly.

Having the relation (\ref{e:Dq-D-Wknown}), we would naturally expect
that it can be inverted, in the sense that
%
\begin{equation}\label{e:D-Dq-Wknown}
F_D = \sum_{m=1}^\infty a_{s,m} F^{*m}_{D_{q,i}|s},
\end{equation}
where $a_s = \{a_{s,m}\}_{m\in \Nset}$ is obtained by reversion of $A_s
= \{A_{s,m}\}_{m\in\Nset}$, that is, their respective
$z$-transformations (or formal power series) $G_{a_s}(z) =
\sum_{m=1}^\infty a_{s,m} z^m$ and $G_{A_s}(z) = \sum_{m=1}^\infty
A_{s,m} z^m$ satisfy $G_{a_s}(G_{A_s}(z)) =G_{A_s}(G_{a_s}(z)) = z$.
The relation (\ref{e:D-Dq-Wknown}) suggests a natural estimator of
$F_D$, defined as (see Section \ref{s:times-sampled-to-renewals})
%
\begin{equation}\label{e:D-Dq-Wknown-estim}
\widehat F_D = \sum_{m=1}^\infty \widehat a_{s,m} \widehat
F^{*m}_{D_{q,i}|s},
\end{equation}
where $\widehat F_{D_{q,i}|s}$ is the empirical distribution function
of $D_{q,i}$ given $W_q=s,$ and $\widehat a_s$ is the reversion of the
sequence $\widehat A_s$, where the latter is defined as $A_s$ by
replacing $f_W$ in its definition by the empirical distribution
$\widehat f_W$. Note that the estimator $\widehat F_D$ is defined for
fixed $s$ and $i$.

\subsection{Contributions to the literature}

Under technical assumptions, we will show that $\widehat F_D$ is an
asymptotically normal estimator of $F_D$, namely,
%
\begin{equation}\label{e:hat-F_D-normal}
    \sqrt{N} (\widehat F_D - F_D) \stackrel{d}{\to} X,
\end{equation}
where $X$ is a limiting Gaussian process and the convergence in
distribution $\stackrel{d}{\to}$ holds in a suitable space of functions
(Theorem \ref{t:an-F_D-sampled}). The approach and some techniques
behind this result are closely related to the works in the so-called
decompounding framework by Buchmann and Gr{\"{u}}bel
\cite{buchmanngrubel2003}, B{\o}gsted and Pitts
\cite{bogstedpitts2009} (see also
\cite{halljuhyun2004,hansenpitts2006}).
These authors consider analogous estimators, but where the sequence
$A_s$ in (\ref{e:Dq-D-Wknown}) is known and hence so is its reversion
sequence $a_s$. We thus deviate from these earlier works by assuming
that $a_s$ also needs to be estimated, which makes the analysis
substantially more complex. More specifically, the limit $X$ in
(\ref{e:hat-F_D-normal}) can be written as
%
\begin{equation}\label{e:hat-F_D-normal-limit}
    X = -\sum_{n=1}^\infty \bigl(a^{(1)}_s*(\zeta\circ a_s)\bigr)_n F_{D_{q,i}|s}^{*n}
    + \Biggl( \sum_{n=1}^\infty a_{s,n} n F_{D_{q,i}|s}^{*(n-1)} \Biggr) * Z,
\end{equation}
where $(\zeta,Z)$ is a Gaussian process, $\circ$ is a composition
operation, $a_s^{(1)}$ is the ``derivative'' of $a_s$ and $*$ is the
convolution. The second term in (\ref{e:hat-F_D-normal-limit}) is the
term found in the decompounding literature when $A_s$ and $a_s$ are
supposed to be known. The first term in (\ref{e:hat-F_D-normal-limit})
is new and accounts  for the variations in $\widehat a_s$ here.

Since $\widehat a_s$ involves the estimator $\widehat f_W$ via
$\widehat A_s$, we also need the asymptotic normality result for
$\widehat f_W$. Although $\widehat f_W$ already appears in the network
traffic literature \cite{hohnveitch2006}, its
statistical properties have not been studied before to the best of our
knowledge, and the asymptotic normality result is also derived here. We
show that, under suitable assumptions,
%
\begin{equation}\label{e:hat-f_W-normal}
    \sqrt{N}(\widehat f_W - f_W) \stackrel{d}{\to} S(\xi),
\end{equation}
where $S(\xi) = \{S(\xi)_w\}_{w\in \Nset}$ is a Gaussian process
(Theorem \ref{t:an-estim-f-W} and Proposition
\ref{p:an-estim-f-W-function}). It is interesting to note that the
result (\ref{e:hat-f_W-normal}) is shown under technical assumptions
which do not cover distributions with heavier tails (such as
heavy-tailed distributions). It is proved, however, that the
assumptions can be thought of as sharp (in the sense of Proposition
\ref{p:necerssary} of Section \ref{s:number-sampled-to-renewals}). What
the asymptotics of $\widehat f_W$ are beyond these assumptions remains
an interesting open question.

We would like to make several other related comments. First, it is
 well known from estimation of flow sizes in the network traffic
context that the performance of the estimator $\widehat f_W$ degrades
rapidly as the sampling probability $q$ decreases. For example, for
$q=0.1$, the inference becomes impractical for reasonable sample sizes
$N$. We shall discuss this fact in light of the derived asymptotic
normality result for $\widehat f_W$ (Sections
\ref{s:number-sampled-to-renewals} and \ref{s:simulations-number}) by
indicating two different regimes for the estimator performance in terms
of the limiting asymptotic variance. The regimes are defined as:
%
\begin{eqnarray}\label{e:2-regimes}
    &\bullet& \mbox{ stable, if }  \sup_{w\in\Nset} ES(\xi)_w^2 <
    \infty;\nonumber\\ [-8pt]\\ [-8pt]
    &\bullet& \mbox{ explosive, if }  \sup_{w\in\Nset}ES(\xi)_w^2 = \infty,\nonumber
  \end{eqnarray}
where $S(\xi)$ is the limiting Gaussian process in
(\ref{e:hat-f_W-normal}) and $ES(\xi)_w^2$ is expressed in terms of
$q$, and $f_W$ or $f_{W_q}$. The performance of the estimator $\widehat
f_W$ is satisfactory in the stable regime, but poor in the explosive
regime. For small $q$, $f_W$ typically belongs to the explosive regime
and hence to the case of poor performance. Analogous difficulties for
smaller $q$ remain when using the estimator~$\widehat F_D$. Because of
these practical considerations, in the network traffic context,
alternative sampling schemes have been  sought and considered, such as
the sample-and-hold method
\cite{cleggetal2008,estanvarghese2002}, and
the dual sampling technique \cite{tuneveitch2008}.
We plan to study inference of interrenewal times in one of these
frameworks and  postpone any real data application to  future work.
This work will therefore be limited to a simulation study (Section
\ref{s:simulation}).

In another direction, the results on characterizing the sampled number
of renewals and sampled interrenewal times, such as the relation
(\ref{e:Dq-D-Wknown}),  contribute to a substantial literature on
sampling of point processes. Sampling (also known as \textit{thinning})
is discussed in several manuscripts such as
\cite{coxisham1980,karr1991,daleyverejones2003}. To the best of our knowledge, sampling
(thinning) of finite renewal processes has been largely unexplored. The
results of our work (Section \ref{s:renewal-to-sampled}) show that
sampling of finite renewal processes does not generally lead to finite
renewal processes (in that, conditionally on the number of sampled
renewals, sampled interrenewal times are generally dependent).

\subsection{Structure of the paper}

The rest of the paper is organized as follows. In Section
\ref{s:preliminaries} we collect the notation used throughout the work
and include other preliminaries. Section \ref{s:renewal-to-sampled}
concerns properties of the finite point process obtained from sampling
the finite renewal process. Inference of the original distributions of
the total number of renewals and interrenewal times is studied in
Section \ref{s:sampled-to-renewal}. A simulation study can be found in
Section \ref{s:simulation}.
 For better readability, most of the proofs are
deferred to Appendix \ref{s:proofs}. The main body of the article
contains the proofs of only those results which we consider key in
individual sections. Finally, Appendix \ref{s:appendix-bounds-taylor}
contains a technical result used in Section
\ref{s:times-sampled-to-renewals}, with bounds on remainder terms in
Taylor expansions of compositions of formal power series.

\section{Notation and other preliminaries}
\label{s:preliminaries}

Here, we introduce notation and make a number of assumptions that will
be used throughout the paper. As in Section \ref{s:intro}, we consider
a finite renewal process consisting of a finite but random number $W$
(with $W\geq 1$) of renewals and positive i.i.d. interrenewal times
$D_i$, $i= 1,\ldots,W-1$, independent of $W$. The terminology of finite
renewal processes here follows  that of
\cite{daleyverejones2003}, Example~5.3(b), page 125. We denote by $D$
the variable with a common distribution of $D_i$. We also let $V$
denote the total duration of the finite renewal process defined by $V =
\sum_{i=1}^{W-1} D_i$ (if $W=1$, then $V=0$). Suppose now that each
renewal, independently of the others, is sampled with a fixed
probability $q\in(0,1)$ to form a sampled finite point process. We let
$W_q$ denote the total number\vspace*{-1pt} of sampled renewals, $D_{q,i}$,
$i=1,\ldots,W_q-1$,  the sampled interrenewal times and $V_q =
\sum_{i=1}^{W_q-1} D_{q,i}$  the total\vspace*{1pt} duration of the sampled finite
point process.

The following notation will be used many times throughout the paper.
For a discrete random variable $X$, its distribution or probability
mass function (p.m.f.) will be denoted by $f_X(x) = P(X=x)$. For
example, we shall write $f_W(w)$, $f_{W_q}(w)$, etc. For a continuous
random variable~$Y$, its distribution function will be denoted\vadjust{\goodbreak} by
$F_Y(y) = P(Y\leq y)$ and its Laplace transform (LT) will be denoted by
$\widetilde F_Y(v) = \int \mathrm{e}^{-vy} F_Y(\mathrm{d}y)$, $v\geq 0$ ($v\in\Rset_+$).
We shall also use the conditional distribution functions $F_{Y|s}(y) =
P(Y\leq y|W_q = s)$ and $F_{Y|s^+}(y) = P(Y\leq y|W_q \geq s)$. So, for
example, we shall write $F_D$, $\widetilde F_D$, $F_{D_{q,i}|s}$, etc.
In several instances, we shall use analogous notation, but where $Y$ is
replaced by a multidimensional.

For a sequence $a = \{a_n\}_{n\in\Nset_0}$, where $\Nset_0 = \Nset \cup
\{0\}$ and $\Nset$ is the set of natural numbers, we denote its formal
power series or its $z$-transformation by $G_a(z) = \sum_{n=0}^\infty
a_n z^n$. Conversely, any such formal power series is associated with a
sequence. When $a$ stands for a~p.m.f. $f_X$, we shall also write $G_X$
in place of $G_{f_X}$. Since a sequence $a$ will not necessarily be
nonnegative, we shall also use the notation $a^+ =
\{a_n^+\}_{n\in\Nset_0}$, defined as $a_n^+ = |a_n|$, and $\G_a(z) =
G_{a^+}(z) = \sum_{n=0}^\infty |a_n|z^n$. In several instances, we
shall use a multidimensional power series $G_a(z_1,\ldots,z_m) =
\sum_{n_1=0}^\infty \cdots \sum_{n_m=0}^\infty a_{{\mathbf n}}
z_1^{n_1}\cdots z_m^{n_m}$ associated with $a = \{a_{\mathbf
n}\}_{{\mathbf n\in \Nset_0^m}}$ and ${\mathbf n} = (n_1,\ldots,n_m)$.

We shall also use the following common operations on sequences (or
their formal power series) $a=\{a_n\}_{n\in\Nset_0}$ and
$b=\{b_n\}_{n\in\Nset_0}$. The composition of $a$ and $b$ will be
denoted $a\circ b$ and is defined by its formal power series as
$G_{a\circ b}(z) = G_a(G_b(z))$. The $k$th derivative of
$a=\{a_n\}_{n\in\Nset_0}$ will be denoted by $a^{(k)} =
\{a^{(k)}_n\}_{n\in\Nset_0}$ and is defined by $G_{a^{(k)}}(z) = \mathrm{d}^k
G_a(z)/\mathrm{d}z^k = \sum_{n=k}^\infty a_n n(n-1)\cdots (n-k+1)z^{n-k}$. The
reversion of a sequence $a = \{a_n\}_{n\in\Nset}$ will be defined as a
sequence $b = \{b_n\}_{n\in\Nset}$ with its formal power series
satisfying $G_b(G_a(z)) = G_a(G_b(z)) = z$. The reversion is defined
for any sequence $a = \{a_n\}_{n\in\Nset}$ with $a_1\neq 0$. As usual,
the symbol $*$ will stand for convolution in either the discrete or
continuous setting. With all these notions for sequences and their
formal power series, we follow the nice monograph of~\cite{henrici1974}.

\section{From finite renewal process to sampled finite point process}
\label{s:renewal-to-sampled}

In this section, we study several characteristics of the sampled finite
point process in terms of the original finite renewal process. We first
briefly consider the number of sampled renewals and then turn to
sampled interrenewal times.

\subsection{Number of sampled renewals}
The relation between the probability mass functions of the number of
sampled renewals and the number of original renewals is given by\vspace*{-1pt}
%
\begin{eqnarray}\label{e:pmf-s_W}
f_{W_q}(s)&=&\sum_{w=s}^\infty P(W_q=s|W=w) P(W=w)\nonumber\\ [-10pt]\\ [-10pt]
 &=& \sum_{w=s}^\infty
{w\choose s} q^{s} (1-q)^{w-s} f_{W}(w),\qquad  s \geq 0.\nonumber
\end{eqnarray}
The function $G_{W_q} (z)$ is given by
\begin{eqnarray}\label{e:rel-G_W:G_s_W}
G_{W_q}(z) &=&\sum_{s=0}^{\infty} z^s \sum_{w=s}^\infty {w\choose s}
q^{s} (1-q)^{w-s} f_{W}(w)\nonumber\\[-2pt]
&=&\sum_{w=1}^{\infty} f_{W}(w) \bigl(zq+(1-q)\bigr)^w \\[-2pt]
&=&G_{W}\bigl(zq+(1-q)\bigr),\qquad |z|<1.\nonumber
\end{eqnarray}
The same relations (\ref{e:pmf-s_W}) and (\ref{e:rel-G_W:G_s_W}) also
appear in \cite{hohnveitch2006} and will be used in
Section \ref{s:number-sampled-to-renewals} to obtain inversion results.

\subsection{Sampled interrenewal times} Given that $s$  renewals have
been sampled, the next  result characterizes the LT of the $i$th $(i<
s)$ sampled interrenewal time  in terms of the  LT of the original
interrenewal times.

\begin{theo}\label{t:LT_sit}
The LT of $D_{q,i}$  given $W_q=s$ can be expressed as
%
\begin{equation}\label{e:LT_sit}
  \widetilde F_{D_{q,i}|s}(v)= G_{A_s} (\widetilde F_D(v)),\qquad  1\leq i < s,  v \geq 0,
 \end{equation}
where $A_s=\{A_{s,m}\}_{m\in \Nset}$ and   $A_{s,m}$ denotes the
probability that the number of original renewals between the $i$th and
$(i+1)$th sampled renewals is equal to $m-1$ given $W_q=s$,
 \begin{equation} \label{eq:LT_sit_coef}
 A_{s,m} = \frac{q^s}{f_{W_q} (s)} \sum_{w=s+m-1}^{\infty} f_W(w)
{w-m\choose s-1}  (1-q)^{w-s}.
 \end{equation}
\end{theo}

\begin{pf}
For $1\leq i < s$ and $t\geq 0$, we have
 \begin{eqnarray}\label{p:LT_sit_aux1}
  \hspace*{-25pt}P (D_{q,i} \leq t, W_{q}=s)&=&\sum_{w=s}^\infty f_W (w) P
(D_{q,i} \leq t, W_{q} = s|W=w) \nonumber \\ [-10pt]\\ [-10pt]
         &=&  \sum_{w=s}^\infty f_W (w) P(W_{q}=s|W=w) P(D_{q,i} \leq t |W=w,W_{q}=s).\nonumber
 \end{eqnarray}
Let $M$ denote the number of original renewals not sampled between the
$i$th and $(i+1)$th sampled renewals plus 1. If the number of original
renewals is $W=w$ and the number of sampled renewals is $W_q=s$, then
$M$ can take the values $1,2,\ldots,w-s+1$, and
\begin{equation}\label{p:LT_sit_aux2}
P(D_{q,i} \leq t |W=w,W_{q}=s) = \sum_{m=1}^{w-s+1} F_{D}^{*m}(t)
P(M=m|W=w,W_{q}=s).
\end{equation}
By considering the location  of the $i$th sampled renewal in the total
number of renewals $w$ along with the possible distinct locations of
the sampled renewals before the $i$th and after the $(i+1)$th sampled
renewals, we obtain that
\begin{eqnarray}\label{p:LT_sit_aux3}
P(M=m|W=w,W_{q}=s) &=&  \sum_{l=i}^{w-m-(s-i-1)}
{l-1\choose i-1}{w-m-l\choose s-i-1}\bigg/{w\choose s} \nonumber \\ [-8pt]\\ [-8pt]
    &=&{w-m\choose s-1}\bigg/{w\choose s},\nonumber
\end{eqnarray}
where $l$ is the location of the  $i$th sampled renewal and the last
equality follows from the identity
%
\begin{equation}\label{eq:ind}
\sum_{l=a}^{x-b} {l-1\choose a-1}{x-l\choose b}={x\choose a+b},\qquad a\geq
1, b \geq 0, a+b \leq x.
\end{equation}
By using (\ref{p:LT_sit_aux1})--(\ref{p:LT_sit_aux3}), we deduce  that
\begin{eqnarray} \label{e:DF_sit}
F_{D_{q,i}|s}(t)&= &\frac{ P (D_{q,i} \leq t, W_{q}=s)} {f_{W_q} (s)}\nonumber\\
& =& \frac{1}{f_{W_q} (s)} \sum_{w=s}^\infty f_W (w) q^{s}
    (1-q)^{w-s} \sum_{m=1}^{w-s+1} {w-m\choose s-1} F^{*m}_D (t)\\
         & = &\sum_{m=1}^{\infty} A_{s,m}  F^{*m}_D (t),\nonumber
\end{eqnarray}
where $A_{s,m}$ is given by (\ref{eq:LT_sit_coef}). The relation
(\ref{e:LT_sit}) follows by taking the LT in (\ref{e:DF_sit}).
\end{pf}

Theorem \ref{t:LT_sit} implies that, when conditioning on $W_q=s$, the
distribution of the $i$th $(i<s)$ sampled interrenewal time depends, in
general, only on $s$ and not on $i$, and hence that the times between
consecutive sampled renewals are identically distributed conditionally
on $W_q$. In contrast to the original finite renewal process which
assumes independence of $W$ and $D_i$, the sampled quantities $W_q$ and
$D_{q,i}$ are dependent.

Other forms of conditioning on the number of sampled renewals are
possible. For example, from (\ref{e:LT_sit}), the LT of $D_{q,i}$
given $W_q\geq s$ can be expressed as
\begin{eqnarray}\label{LT:geq-s}
 \widetilde F_{D_{q,i}|s^+}(v) &=& \sum_{s^{\prime}=s}^{\infty}  \widetilde F_{D_{q,i}|s^{\prime}}(v)
P(W_q=s^{\prime}|W_q \geq s) \nonumber \\
&=&\frac{1}{P(W_q \geq s)}\sum_{s^{\prime}=s}^{\infty}
\sum_{m=1}^{\infty} A_{s^{\prime},m}
\widetilde F_{D} (v)^m f_{W_q}(s^{\prime})\\
&=& G_{A_{s^+}} (\widetilde F_{D} (v)),\nonumber
\end{eqnarray}
where $A_{s^+}=\{A_{s^+,m}\}_{m \in \Nset}$ and $A_{s^{+},m}$ is the
probability that the number of original renewals between the $i$th and
$(i+1)$th sampled renewals is equal to $m-1$ given $W_q\geq s$,
\begin{eqnarray}\label{series-coef:geq-s}
 \hspace*{-25pt}  A_{s^{+},m} &=& \frac{1}{P(W_q \geq s)} \sum_{s^{\prime}=s}^{\infty}
\sum_{w=s^{\prime}+m-1}^{\infty} f_W(w) {w-m\choose s^{\prime}-1}
q^{s^{\prime}} (1-q)^{w-s^{\prime}} \nonumber \\ [-8pt]\\ [-8pt]
  &=& \frac{q (1-q)^{m-1}}{P(W_q \geq s)} \sum_{w=m+s-1}^\infty f_W(w)
 \Biggl( 1- \sum_{s^{\prime}=0}^{s-2}  {w-m \choose s^{\prime}}
q^{s^{\prime}} (1-q)^{w-m-s^{\prime}} \Biggr).\nonumber
\end{eqnarray}
In the rest of the paper, we shall focus only on the conditioning
$W_q=s,$ but similar results can be derived for other forms of
conditioning such as $W_q \geq s$.

The next result gives the joint LT for a finite number of sampled
interrenewal times. We shall use this general result below to
investigate dependence between sampled interrenewal times (see Section
\ref{s:examples}); see Appendix \ref{s:proofs} for the proof.

\begin{theo}\label{t:LT_joint_sit}
The joint LT of $\mathbf D_{q,n}=(D_{q,i_1},\ldots,D_{q,i_n})$ with $1
\leq i_1<\cdots< i_n$ given $W_q=s$, $s >i_n$, can be expressed as
\begin{equation} \label{e:LT_joint_sit-1}
\widetilde F_{\mathbf D_{q,n}|s}(\mathbf v) = G_{B_s} (\widetilde
F_D(v_1),\ldots,\widetilde F_D(v_n)),\qquad  \mathbf v = (v_1,\ldots,v_n) \in
\Rset^n_+,
\end{equation}
where  $B_s=\{B_{s,\mathbf m}\}_{\mathbf m \in \Nset^n}$ with $\mathbf
m = (m_1,\ldots,m_n)$ and $B_{s,\mathbf m}$ denotes the probability
that the number of renewals between the $i_j$th and $(i_{j}+1)$th
sampled renewals is equal to $m_j-1$, $1 \leq j \leq n$, given $W_q=s$:
\begin{equation}\label{e:LT_joint_sit_coef}
  B_{s,\mathbf m} = \frac{q^s}{f_{W_q} (s)} \sum_{w= s+m-n}^\infty f_W(w)
{w- m\choose s-n} (1-q)^{w-s}
\end{equation}
with $m = m_1+\cdots + m_n$.
\end{theo}

By Theorem \ref{t:LT_joint_sit}, the joint distribution of a subset of
sampled interrenewal times given $W_q=s$ depends only on $s$ and not on
the indices  of the sampled times considered. Another relation of
interest is between the duration of the sampled finite point process
and the original interrenewal times; see Appendix \ref{s:proofs} for
the proof.

\begin{prop} \label{p:LT_s_V}
The LT of  $V_q$ given that $W_q\geq 2$ can be expressed as
%
\begin{equation}\label{e:LT_s_V}
  \widetilde F_{V_q|2^+}(v)= G_C (\widetilde F_D(v)),\qquad v\geq 0,
 \end{equation}
where  $C=\{C_m\}_{m\in \Nset}$ and $C_m$ is the probability that the
number of renewals between the first and last sampled renewals is equal
to $m-1$, given $W_q\geq 2$:
%
\begin{equation}\label{e:LT_sit_coef}
  C_{m} = \frac{q^2}{P(W_q \geq 2)} \sum_{w=m+1}^\infty f_W(w) (w-m) (1-q)^{w-m-1}.
\end{equation}
\end{prop}


Since the LT of the duration of the original finite renewal process
$V=\sum_{i=1}^{W-1} D_i$, given\break $W\geq 2$, can be expressed as
\[
\widetilde F_{V|W\geq 2} (v) = \frac{1}{P(W\geq 2)} \sum_{w=1}^{\infty}
f_W (w+1) \widetilde F_D (v)^w,
\]
we expect that
%
\begin{equation}\label{LT:it}
\widetilde F_D (v) = \sum_{n=1}^{\infty} D_n \widetilde F_{V|W\geq
2}(v)^n,
\end{equation}
where $D=\{D_n\}_{n\in\Nset}$ is the reversion of $\{f_W(w+1)/P(W\geq
2)\}_{w\in\Nset}$ (defined in Section \ref{s:preliminaries}). By
plugging $(\ref{LT:it})$ into (\ref{e:LT_s_V}), we can obtain an
expression for the duration $V_q$ of the sampled finite point process
in terms of the duration $V$ of the finite renewal process.

Also, note that, equivalently, the relations (\ref{e:LT_sit}),
(\ref{e:LT_joint_sit-1}) and (\ref{e:LT_s_V}) can be expressed in terms
of the distributions functions (see (\ref{e:DF_sit}), and relations
(\ref{e:DF_joint_sit}) and (\ref{e:DF_s_V}) in Appendix \ref{s:proofs},
respectively) or via the characteristic functions.

\subsubsection{Examples}
\label{s:examples} In the following, we examine characteristics of the
sampled finite point process for particular distributions of $W$.

\subsubsection*{Geometric distribution}
Suppose that  $W$ is a geometric random variable with parameter $c$,
that is, $f_W(w)=c^{w-1}(1-c)$, $w\geq 1$, $c \in (0,1)$. Substituting
this $f_W(w)$ into (\ref{eq:LT_sit_coef}), we  obtain, after
straightforward calculations, that
\[
  A_{s,m}=  \bigl(c(1-q)\bigr)^{m-1} \bigl(1-c(1-q)\bigr),\qquad  m\geq 1.
\]
Therefore,  $A_s$ is the p.m.f. of a geometric distribution with
parameter $c(1-q)$ which does not depend on $s$. From (\ref{e:LT_sit})
and using the generating function of the geometric distribution, we
have
\[
\widetilde F_{D_{q,i}|s}(v)=\frac{(1-c(1-q)) \widetilde F_D
(v)}{1-c(1-q) \widetilde F_D(v)}.
\]
Now, substituting $f_W(w)$ into (\ref{e:LT_joint_sit_coef}), we
conclude, after some algebra, that
\[
 B_{s,\mathbf m} = A_{s,m_1}\cdots A_{s,m_n},\qquad  \mathbf m = (m_1,\ldots,m_n)\in
 \Nset^n,
\]
and hence, from (\ref{e:LT_joint_sit-1}),
\[
  \widetilde  F_{\mathbf D_{q,n}|s} (\mathbf v)=\prod_{i=1}^n \frac{(1-c(1-q))
\widetilde F_D(v_i)}{1-c(1-q) \widetilde F_D (v_i)}.
\]
The latter expression shows that  when $W$ is geometrically
distributed, the sampled interrenewal times are  independent and
identically distributed (i.e., the sampled finite point process is a
finite renewal process). The LT of $V_q$ in (\ref{e:LT_s_V}), where $C$
is now the p.m.f. of a geometric distribution with parameter $1-c$,
simplifies to
\[
 \widetilde F_{V_q|2^+}(v)=\frac{ c \widetilde F_D(v)}{1-(1-c) \widetilde F_D (v)}.
\]

\subsubsection*{Pareto distribution}
Suppose now that $W$ is a discrete Pareto distribution, that is,
$f_W(w)=w^{-\alpha-1}/\zeta(\alpha+1)$, $w\geq 1$, $\alpha>0$ and
$\zeta(z)=\sum_{w=1}^\infty w^{-z}$ is the Riemann zeta function.
Unlike in the previous example, closed forms for LT's of the sampled
interrenewal times are not available. Nevertheless, we can show that
for a particular example, the sampled interrenewal times are not
conditionally independent. Taking $s=3$, $i_1=1$ and $i_2=2$, the
random variables $D_{q,1}$ and $D_{q,2}$ are independent given $W_q=3$
if and only if
%
\begin{equation}
G_{B_3} (\widetilde F_D(v_1),\widetilde F_D(v_2)) =G_{A_3}(\widetilde
F_D(v_1)) G_{A_3}(\widetilde F_D(v_2)) \label{ind:pareto}
\end{equation}
for all $(y_1,y_2) \in \Rset_+^2$. Since $\widetilde F_D(y)$ is
continuous in $y$, it takes all the values in $(0,1)$. Then,
(\ref{ind:pareto}) implies that $G_{B_3} (w_1,w_2) =G_{A_3}(w_1)
G_{A_3}(w_2)$ for all $w_1,w_2\in(0,1)$ and\vspace*{1pt} hence $B_{3,\mathbf
m}=A_{3,m_1}A_{3,m_2}$ for all $\mathbf m = (m_1,m_2) \in \Nset^2$. The
coefficients of the term $\widetilde F_D(v_1)\widetilde F_D(v_2)$ on
the right- and left-hand sides of (\ref{ind:pareto}) are, respectively,
\[
{A_{3,1}^2}=
\frac{9(\operatorname{Li}_{\alpha-1}(1-q)-3\operatorname{Li}_{\alpha}(1-q)+2\operatorname{Li}_{\alpha+1}(1-q))^2}
{\operatorname{Li}_{\alpha-2}(1-q) - 3\operatorname{Li}_{\alpha-1}(1-q)+ 2
\operatorname{Li}_{\alpha}(1-q)}
\]
and
\[
 B_{3,(1,1)}= \frac{6(1 - q + \operatorname{Li}_{\alpha}(1-q) -2 \operatorname{Li}_{\alpha+1}(1-q))}
{\operatorname{Li}_{\alpha-2}(1-q) - 3\operatorname{Li}_{\alpha-1}(1-q)+ 2
\operatorname{Li}_{\alpha}(1-q)},
\]
where $\operatorname{Li}_{n}(a)=\sum_{w=1}^{\infty} a^w/w^n$, $|a|<1$, is the
polylogarithm function. Since ${A_{3,1}^2} \neq B_{3,(1,1)}$, the
random variables $D_{q,1}$ and $D_{q,2}$ are not conditionally
independent. This is also why we use the term ``sampled finite point
process'' throughout the paper instead of ``sampled finite renewal
process''.

\subsubsection*{Heavy-tailed distributions}

 Consider the case of a heavy-tailed
distribution for $W,$ in the sense that
%
\begin{equation}\label{e:ht}
    f_W(w) \sim c \alpha w^{-\alpha-1}
\end{equation}
as $w\to \infty$, where $\alpha\in (1,2)$ and $c>0$. Such distributions
are common models for long flow sizes in network traffic studies (e.g.,
\cite{chabchoubfrickerguilleminrobert2010}). As in
the example with Pareto distribution above, closed forms for LT's of
the sampled interrenewal times are not available under (\ref{e:ht}). We
can nevertheless provide a number of interesting qualitative results
concerning the case (\ref{e:ht}).

The relation (\ref{e:ht}) implies that
%
\begin{equation}\label{e:ht2}
   P(W >w)  \sim c w^{-\alpha}
\end{equation}
as $w\to \infty$. As indicated in \cite{hohnveitch2006}, page
70, by the results of \cite{binghamgoldieteugels1987}, page
333, this implies that
%
\begin{equation}\label{e:ht2-q}
    P(W_q >w) \sim q^\alpha c w^{-\alpha},
\end{equation}
that is, $W_q$ is also heavy-tailed with the same parameter $\alpha$.
(Note that (\ref{e:ht2-q}) does not, in general, imply that
$f_{W_q}(w)\sim q^\alpha c \alpha w^{-\alpha-1}$.) Also, note that in
the case (\ref{e:ht2}), the total duration $V$ is expected to be
heavy-tailed, even for light-tailed interrenewal times $D_i$. Indeed,
by Robert and Segers \cite{robertsegers2008}, Theorem 3.2, if
(\ref{e:ht2}) holds and
\begin{equation}\label{e:ht-cond}
    P(D_1>t) = \mathrm{o}(t^{-\alpha}),
\end{equation}
then
%
\begin{equation}\label{e:ht-duration}
    P(V>t) \sim  P  \biggl(W> \frac{t}{ED_1} \biggr) \sim
 c (ED_1)^\alpha t^{-\alpha},
\end{equation}
that is, the tail of total duration is dominated by that of $W$.

On the sampling side, the total duration is also characterized by heavy
tails. Indeed, by Proposition \ref{p:LT_s_V} above and \cite{robertsegers2008}, Theorem 3.2,
%
\begin{equation}\label{e:ht-duration-q}
    P(V_q>t) \sim P  \biggl(C > \frac{t}{ED_1} \biggr) \sim
    c_0 (ED_1)^\alpha t^{-\alpha},
\end{equation}
where a random variable $C$  is such that $P(C=m) = C_m$, where $C_m$
is defined in (\ref{e:LT_sit_coef}), as long as (\ref{e:ht-cond}) holds
and
\begin{equation}\label{e:ht-C}
    C_m \sim c_0 \alpha m^{-\alpha-1}.
\end{equation}
To show (\ref{e:ht-C}), observe that for large $m$, using (\ref{e:ht}),
\begin{eqnarray}\label{e:ht-C-proof}
C_m &\sim& \frac{q^2 c\alpha }{P(W_q \geq 2)} \sum_{w=m+1}^\infty w^{-\alpha-1}
  (w-m) (1-q)^{w-m-1} \nonumber \\
   &=&  m^{-\alpha-1} \frac{q^2 c\alpha }{P(W_q \geq 2)} \sum_{k=1}^\infty \biggl(1 +
\frac{k}{m} \biggr)^{-\alpha-1}
  k (1-q)^{k-1} \nonumber \\
   &\sim& m^{-\alpha-1} \frac{q^2 c\alpha }{P(W_q \geq 2)} \sum_{k=1}^\infty k
   (1-q)^{k-1}\\
   &=& m^{-\alpha-1} \frac{q^2 c\alpha }{P(W_q \geq 2)} \Biggl(- \frac{\mathrm{d}}{\mathrm{d}q}
\sum_{k=0}^\infty
   (1-q)^k \Biggr)\nonumber\\
   &=& m^{-\alpha-1} \frac{q^2 c\alpha }{P(W_q \geq 2)}\biggl (- \frac{\mathrm{d}}{\mathrm{d}q} \frac{1}{q}
\biggr)
   = m^{-\alpha-1} \frac{c\alpha }{P(W_q \geq 2)} =:
    m^{-\alpha-1} c_0\alpha.\nonumber
\end{eqnarray}

Note also that the argument (\ref{e:ht-C-proof}) above does not apply
to the sampled interrenewal times. For example, for the first sampled
interrenewal time $D_{q,1}$, the analogous argument would show that,
for coefficients in (\ref{eq:LT_sit_coef}),
\begin{equation}\label{e:ht-not}
    A_{s,m} \sim c (1-q)^m m^{-\alpha-1}
\end{equation}
as $m\to \infty$, which is not heavy-tailed, in contrast to
(\ref{e:ht-C}).


\section{From sampled finite point process to finite renewal process}
\label{s:sampled-to-renewal}

In this section, we study inference of the original distributions of
the number of renewals and interrenewal times.

\subsection{Distribution of number of renewals}
\label{s:number-sampled-to-renewals}

We are interested here in estimating the p.m.f. $f_W(w)$ from i.i.d.\
observations of $W_q$. We revisit a nonparametric estimator of $f_W(w)$
introduced in \cite{hohnveitch2006} and clarify
several issues surrounding its use and properties. A number of open
questions will also be raised.

Estimation of $f_W(w)$ is based on a theoretical inversion of the
relation (\ref{e:pmf-s_W}) which will be discussed first. The relation
(\ref{e:rel-G_W:G_s_W}) can be written as $G_W(z) = G_{W_q}(q^{-1}z +
(1-q^{-1})),$ which has the same form as the original relation
(\ref{e:rel-G_W:G_s_W}) when  $W$ and $W_q$ are exchanged, and $q$ is
replaced by $q^{-1}$. In view of (\ref{e:pmf-s_W}), we would then
expect that
\begin{eqnarray}\label{e:f_Wq->f_W}
   f_W(w) &=& \sum_{s=w}^\infty {s\choose w} (q^{-1})^w (1-q^{-1})^{s-w} f_{W_q}(s) \nonumber
   \\ [-8pt]\\ [-8pt]
   &=& \sum_{s=w}^\infty {s\choose w} \frac{(-1)^{s-w}}{q^s} (1-q)^{s-w}
   f_{W_q}(s),\qquad w\geq 1.\nonumber
\end{eqnarray}
Hohn and Veitch \cite{hohnveitch2006} claim that (\ref{e:f_Wq->f_W})
holds when $q\in (0.5,1)$ (with this choice, note that $q^{-s}(1-q)^s
\in (0,1)$ in (\ref{e:f_Wq->f_W})). As the following elementary result
shows, this is not a necessary condition. The relation
(\ref{e:f_Wq->f_W}) also holds for $q\in (0,0.5]$ as long as
$f_{W_q}(w)$ or $f_W(w)$ decays to zero fast enough; see Appendix
\ref{s:proofs} for the proof.

\begin{prop}\label{p:f_Wq->f_W-1}
If
\begin{equation}\label{e:f_Wq->f_W-1-cond-main}
\sum_{s=n}^\infty {s\choose n} \frac{(1-q)^{s-n}}{q^s} f_{W_q}(s) =
\sum_{w=n}^\infty {w\choose n} 2^{w-n} (1-q)^{w-n} f_W(w) < \infty,\qquad
n\geq 1,
\end{equation}
then the relation (\ref{e:f_Wq->f_W}) holds.
\end{prop}


For example, for geometric $W$ satisfying $f_W(w) = c^{w-1}(1-c)$,
$w\geq 1$, $c\in(0,1)$, the condition (\ref{e:f_Wq->f_W-1-cond-main})
holds when $2-(2-\epsilon)q<c^{-1}$, where $\epsilon >0$ is arbitrarily
small and where we have used the fact that ${w\choose n}$ in
(\ref{e:f_Wq->f_W-1-cond-main}) can be bounded by $w^n$ up to a
multiplicative constant. Note also that (\ref{e:f_Wq->f_W-1-cond-main})
always holds for $q\in(0.5,1)$. When $q\in (0,0.5]$, on the other hand,
a number of distributions $f_W(w)$ of interest, such as heavy-tailed
distributions (Section \ref{s:examples}), do not satisfy
(\ref{e:f_Wq->f_W-1-cond-main}). In these cases, and generally for
$q\in (0,0.5]$, the relation (\ref{e:f_Wq->f_W}) needs to be modified
by a procedure used and referred to as an analytic continuation in \cite{hohnveitch2006}, page 71. The procedure is defined
below. We should note that this procedure is viewed below  as a
convenient algebraic trick that makes series converge (see, e.g., the
proof of Proposition \ref{p:f_Wq->f_W-2}) rather than as a suitable
analytic continuation, which is a complementary viewpoint followed in
\cite{hohnveitch2006}.

Let $z_0 = 1-q$ and pick an arbitrary sequence $z_k$, $k=1,\ldots,l$,
such that $1>z_0>z_1>\cdots >z_{l-1}>z_l=0$ and $z_k\in C_{k-1} =
\{z\dvtx|z-z_{k-1}| < 1 - z_{k-1}\}$, $k=1,\ldots,l$. For a sequence $x =
\{x_n\}_{n\in\Nset}$, define formally and recursively sequences
$T^{(k)}(x) = \{T^{(k)}(x)_n\}_{n\in\Nset}$, $k=1,\ldots,l$, as
%
\begin{equation}\label{e:T^k_x}
    T^{(k)}(x)_n = \sum_{i=n}^\infty {i\choose n} T^{(k-1)}(x)_i (z_k -
    z_{k-1})^{i-n},\qquad n\geq 1,
\end{equation}
where $T^{(0)}(x)_i = x_i/q^i$. It is also convenient to define the
mapping underlying (\ref{e:f_Wq->f_W}), that is, $S(x) =
\{S(x)_n\}_{n\in\Nset}$, where
%
\begin{equation}\label{e:S_x}
    S(x)_n = \sum_{i=n}^\infty {i\choose n} \frac{(-1)^{i-n}}{q^i}
    (1-q)^{i-n} x_i,\qquad  n\geq 1.
\end{equation}
The following elementary result relates $S(x)$ and $T^{(l)}(x)$ when
$x$ satisfies a natural condition.

\begin{prop}\label{p:T-S}
If a sequence $x = \{x_n\}_{n\in\Nset}$ satisfies
\begin{equation}\label{e:T-S-cond}
    \sum_{i=n}^\infty {i\choose n} \frac{(1-q)^{i-n}}{q^i} |x_i| <
    \infty,\qquad n\geq 1,
\end{equation}
then
\begin{equation}\label{e:T=S}
    T^{(l)}(x) = S(x),
\end{equation}
where $T^{(l)}(x)$ and $S(x)$ are defined in (\ref{e:T^k_x}) and
(\ref{e:S_x}), respectively.
\end{prop}


The next result formalizes the fact that $f_W$ can be obtained as
$T^{(l)}(f_{W_q})$. Note that this is true without any assumptions on
$q$ and $f_W$.

\begin{prop}\label{p:f_Wq->f_W-2}
For any $q\in(0,1)$ and p.m.f. $f_W$, we
have\begin{equation}\label{e:f_Wq->f_W-2} f_W = T^{(l)}(f_{W_q}),
\end{equation}
where $T^{(l)}$ is defined in (\ref{e:T^k_x}).
\end{prop}


With theoretical inversion formulas (\ref{e:f_Wq->f_W}) and
(\ref{e:f_Wq->f_W-2}), we can now turn to estimation. Let $W_{q,k}$,
$k=1,\ldots,N$, be i.i.d. copies of the variable $W_q$ and
%
\begin{equation}\label{e:estim-f_Wq}
    \widehat f_{W_q}(s) = \frac{1}{N} \sum_{k=1}^N 1_{\{W_{q,k} = s\}},\qquad s\geq 0,
\end{equation}
be the empirical p.m.f. of $f_{W_q}$, where $1_A$ denotes an indicator
function of an event $A$. Note that we assume, in particular, that the
event $W_{q,k}=0$ can be observed in practice (or, equivalently, $N$~and $W_{q,k}\geq 1$ can be observed in practice). In the network
traffic context, we naturally observe only $W_{q,k}\geq 1$. The total
number of flows $N$ is deduced from the additional information in the
sampled packet headers (e.g.,
\cite{duffieldlundthorip2002}). (To be more precise, $N$ is actually
estimated, but we suppose it to be known for the sake of simplicity.)

In view of (\ref{e:f_Wq->f_W-2}), it is natural to introduce the
following nonparametric estimator of $f_W$:
\begin{eqnarray}\label{e:estim-f_W-1}
  \widehat f_W(w) &=& T^{(l)}(\widehat f_{W_q})_w \\ \label{e:estim-f_W-2}
   &=& S(\widehat f_{W_q})_w = \sum_{s=w}^\infty {s\choose w} \frac{(-1)^{s-w}}{q^s}
   (1-q)^{s-w} \widehat f_{W_q}(s),\qquad
    w\geq 1,
\end{eqnarray}
where $T^{(l)}(x)$ and $S(x)$ are defined in (\ref{e:T^k_x}) and
(\ref{e:S_x}), respectively. The first equality in
(\ref{e:estim-f_W-2}) follows from Proposition \ref{p:T-S} since only a
finite number of the  $\widehat f_{W_q}(s)$'s are nonzero.

We will next show that the estimator (\ref{e:estim-f_W-1}) is
asymptotically normal under suitable assumptions. The suitable
assumptions are quite strong, but we do not expect that they can be
weakened much, as explained in Proposition \ref{p:necerssary} below and
a discussion surrounding it. For $w\geq 1$, we also let
\begin{eqnarray}\label{e:R_w-1}
  R_{q,w} &=&  \sum_{s=w}^\infty {s\choose w}^2 \frac{(1-q)^{2(s-w)}}{q^{2s}} f_{W_q}(s) \\\label{e:R_w-2}
   &=& \sum_{k=w}^\infty f_W(k) (1-q)^{k-2w} {k\choose w} \sum_{s=w}^k {s\choose w} {k-w\choose s-w} (q^{-1}-1)^s,
\end{eqnarray}
where the second equality follows from (\ref{e:pmf-s_W}).

\begin{theo}\label{t:an-estim-f-W}
Suppose that
%
\begin{equation}\label{e:an-estim-f-W-cond-main}
     R_{q,w} < \infty,\qquad  w\geq 1,
\end{equation}
where $R_{q,w}$ is defined in (\ref{e:R_w-1}). Then, as $N\to\infty$,
%
\begin{equation}\label{e:an-estim-f-W}
    \bigl\{\sqrt{N} \bigl(\widehat f_W(w) - f_W(w) \bigr)  \bigr\}_{w\in\Nset}
    \rightarrow \{ S(\xi)_w \}_{w\in\Nset},
\end{equation}
where the convergence is in the sense of  finite-dimensional
distributions, $S(x)$ is defined in (\ref{e:S_x}) and $\xi =
\{\xi_s\}_{s\in\Nset}$ is a zero-mean Gaussian process with the
covariance structure
%
\begin{equation}\label{e:xi-cov-struct}
    E\xi_{s_1}\xi_{s_2} = f_{W_q}(s_1) \delta_{s_1,s_2} - f_{W_q}(s_1)
    f_{W_q}(s_2)
\end{equation}
(and, as usual, $\delta_{s_1,s_2} = 1$ if $s_1 = s_2$, and $=0$ if
$s_1\neq s_2$, this being the Kronecker symbol). In particular, the
limiting variables $S(\xi)_w$ are zero-mean, Gaussian and have the
variance
\begin{equation}\label{e:an-var-limit}
    E S(\xi)_w^2 = R_{q,w} - f_W(w)^2.
\end{equation}
\end{theo}

\begin{pf}
We only consider the convergence (\ref{e:an-estim-f-W}) at fixed $w\geq
1$. Note that, using Proposition \ref{p:T-S}, $\sqrt{N} (\widehat
f_W(w) - f_W(w) ) = S(\sqrt{N} (\widehat f_{W_q} - f_{W_q}))_w$. For
fixed $j\geq w$ and $x = \{x_n\}_{n\in\Nset}$, define
\[
S(x)_{j,w} = \sum_{i=w}^j {i\choose w} \frac{(-1)^{i-w}}{q^i}
(1-q)^{i-w} x_i.
\]
Using \cite{billingsley1999}, Theorem 3.2, page 28, it is
enough to show that:
\begin{longlist}
    \item  $S(\sqrt{N} (\widehat f_{W_q} - f_{W_q}))_{j,w} \stackrel{d}{\to}
    S(\xi)_{j,w}$ as $N\to \infty$;
    \item $S(\xi)_{j,w} \stackrel{d}{\to} S(\xi)_{w}$ as $j\to \infty$;
    \item for any $\delta >0$,
\[
\limsup_{j\to \infty} \limsup_{N\to \infty} P\bigl( \bigl|S\bigl(\sqrt{N} (\widehat
f_{W_q} - f_{W_q})\bigr)_{j,w} - S\bigl(\sqrt{N} (\widehat f_{W_q} -
f_{W_q})\bigr)_{w}\bigr| > \delta\bigr) = 0.
\]
\end{longlist}
The convergence in (i) is elementary since $\{\sqrt{N} (\widehat
f_{W_q}(s) - f_{W_q}(s)) \}_{s\in\Nset}$ converges to $\xi =
\{\xi_s\}_{s\in\Nset}$ in the sense of finite-dimensional distributions
and since, for fixed $j$, $S(x)_{j,w}$ involves only a~finite number of
elements of $x=\{x_n\}_{n\in\Nset}$.

The convergence in (ii) can be proven in a stronger sense, that of
almost sure convergence. For this, observe that in the sense of
finite-dimensional distributions,
%
\begin{eqnarray}\label{e:xi-repres}
  \{\xi_s\}_{s\in\Nset}  &\stackrel{d}{=} &  \bigl\{B(F_{W_q}(s)) - B\bigl(F_{W_q}(s-1)\bigr)
    \bigr\}_{s\in\Nset} + \{f_{W_q}(s)B(1)\}_{s\in\Nset} \nonumber \\ [-8pt]\\ [-8pt]
   &=:&\bigl\{\xi^{(1)}_s\bigr\}_{s\in\Nset} +
    \bigl\{\xi^{(2)}_s\bigr\}_{s\in\Nset},\nonumber
\end{eqnarray}
where $B = \{B(t)\}_{t\in [0,1]}$ is a standard Brownian motion. It is
then enough to show almost sure convergence of the corresponding terms
for $\xi^{(1)}$ and $\xi^{(2)}$. Doing this for $\xi^{(2)}$ is
elementary and thus we do so only for $\xi^{(1)}$. Note that
\begin{equation}\label{e:xi-1-repres}
 \bigl\{\xi^{(1)}_s\bigr\}_{s\in\Nset}  \stackrel{d}{=}
 \{f_{W_q}(s)^{1/2} \eta_s\}_{s\in\Nset},
\end{equation}
where $\eta_s$ are i.i.d. $N(0,1)$ random variables. By the three
series theorem, $S(\xi^{(1)})_{j,w} \to S(\xi^{(1)})_w$ a.s. as long as
the condition (\ref{e:an-estim-f-W-cond-main}) holds.

For the convergence in (iii), note that the probability in (iii)
 can be expressed and bounded as
\begin{eqnarray*}
&&P\Biggl( \Biggl| \sum_{s=j+1}^\infty {s\choose w} \frac{(-1)^{s-w}}{q^s}
(1-q)^{s-w} \sqrt{N} \bigl(\widehat f_{W_q}(s) - f_{W_q}(s)\bigr) \Biggr|
> \delta \Biggr)\\
&&\quad\leq \delta^{-2} \sum_{s=j+1}^\infty {{s\choose w}}^2
\frac{(1-q)^{2(s-w)}}{q^{2s}}
 E\bigl(\sqrt{N} \bigl(\widehat f_{W_q}(s) - f_{W_q}(s)\bigr)\bigr)^2\\
&&\qquad{}+ 2 \delta^{-2} \sum_{j+1\leq s_1 <s_2} {s_1\choose w} {s_2\choose w}
\frac{(1-q)^{s_1+s_2-2w}}{q^{s_1+s_2}}\\
&&\qquad\hphantom{{}+ 2 \delta^{-2} \sum_{j+1\leq s_1 <s_2}}{}\times\bigl| E N \bigl(\widehat f_{W_q}(s_1) - f_{W_q}(s_1)\bigr) \bigl(\widehat f_{W_q}(s_2) -
f_{W_q}(s_2)\bigr) \bigr|\\
&&\quad= \delta^{-2} \sum_{s=j+1}^\infty {{s\choose w}}^2
\frac{(1-q)^{2(s-w)}}{q^{2s}}
 \bigl( f_{W_q}(s) - f_{W_q}(s)^2 \bigr)\\
&&\qquad{}+ 2 \delta^{-2} \sum_{j+1\leq s_1 <s_2} {s_1\choose w} {s_2\choose w}
\frac{(1-q)^{s_1+s_2-2w}}{q^{s_1+s_2}}
f_{W_q}(s_1) f_{W_q}(s_2) \leq \delta^{-2} (Q_j + 2 Q_j^2),
\end{eqnarray*}
where $Q_j =  \sum_{s=j+1}^\infty {{s\choose w}}^2
\frac{(1-q)^{2(s-w)}}{q^{2s}} f_{W_q}(s)$. It remains to observe that
$Q_j\to 0$ as $j\to \infty$, by the assumption
(\ref{e:an-estim-f-W-cond-main}). The proof of (\ref{e:an-var-limit})
is now elementary.
\end{pf}

Note that the condition (\ref{e:an-estim-f-W-cond-main}) is quite
strong and, for example, does not allow for heavy-tailed distributions
when $q\in (0,0.5]$. When this condition does not hold, since $\sqrt{N}
(\widehat f_W(w) - f_W(w) ) = T^{(l)}(\sqrt{N} (\widehat f_{W_q} -
f_{W_q}))_w$, we may think that the limit of the left-hand side of
(\ref{e:an-estim-f-W}) should be
%
\begin{equation}\label{e:an-estim-f-W-expected-limit}
    T^{(l)}(\xi) = \bigl\{ T^{(l)}(\xi)_w \bigr\}_{w\in\Nset},
\end{equation}
where the process $\xi$ is as in Theorem \ref{t:an-estim-f-W}. This,
however, is not expected to hold. In fact, without the condition
(\ref{e:an-estim-f-W-cond-main}), the expected limit
(\ref{e:an-estim-f-W-expected-limit}) is not even well defined, as the
result below shows. Consider the case of $l=2$, $z_2=0$ and
$T^{(2)}(\xi)$ for simplicity (the argument can also be extended to
general $l$). With the representation (\ref{e:xi-repres}), the term
$T^{(2)}(\xi^{(2)})$ is well defined by arguing as in the proof of
Proposition \ref{p:f_Wq->f_W-2}. We will show that without the
condition (\ref{e:an-estim-f-W-cond-main}), $T^{(2)}(\xi^{(1)})$ is not
well defined, in the following natural sense. Note that
\[
T^{(1)}\bigl(\xi^{(1)}\bigr)_n = \sum_{s=n}^\infty {s\choose n}
\frac{(z_1-(1-q))^{s-n}}{q^s} f_{W_q}(s)^{1/2} \eta_s
\]
is well defined, where we take $\xi^{(1)}_s = f_{W_q}(s)^{1/2} \eta_s$
with $\eta_s$ as in (\ref{e:xi-1-repres}) and let
\[
T^{(2)}\bigl(\xi^{(1)}\bigr)_{j,w} = \sum_{k=w}^j {k\choose w} (-z_1)^{k-w}
T^{(1)}\bigl(\xi^{(1)}\bigr)_k.
\]
The proof of the following result can be found in Appendix
\ref{s:proofs}.

\begin{prop}\label{p:necerssary}
The condition (\ref{e:an-estim-f-W-cond-main}) is necessary for
$T^{(2)}(\xi^{(1)})_{j,w}$ to have a limit in distribution as
$j\to\infty$.
\end{prop}

%

We shall conclude this section by considering the convergence
(\ref{e:an-estim-f-W}) in a suitable space of sequences. We let
$l_{\infty,a} = \{x = \{x_n\}_{n\in\Nset_0}\dvtx \|x\|_{\infty,a} : =
\sup_{n\geq 0} a^n |x_n|<\infty \}$. See Appendix \ref{s:proofs} for
the proof.

\begin{prop}\label{p:an-estim-f-W-function}
If, for $b'>0$,
%
\begin{equation}\label{e:an-estim-f-W-function-cond}
    \sum_{s=0}^\infty f_{W_q}(s)^{1/2} \biggl( \frac{b'+1-q}{q} \biggr)^s < \infty,
\end{equation}
then the convergence (\ref{e:an-estim-f-W}) also holds in the space
$l_{\infty,b}$ for any $0<b<b'$.
\end{prop}

Note that the condition (\ref{e:an-estim-f-W-function-cond}) is
stronger than (\ref{e:an-estim-f-W-cond-main}).

%
%

Theorem \ref{t:an-estim-f-W} and Proposition
\ref{p:an-estim-f-W-function} suggest several possibilities for
confidence intervals (CI's) of $f_W$. For example, following Theorem
\ref{t:an-estim-f-W}, we can define the 100$\alpha\%$ CI of $f_W(w)$
for fixed $w$ as
%
\begin{equation}\label{e:ci-1}
    \bigl(\widehat f_W(w) - z_\alpha N^{-1/2} (\widehat{ES(\xi)_w^2})^{1/2},
    \widehat f_W(w) + z_\alpha N^{-1/2} (\widehat{ES(\xi)_w^2})^{1/2}\bigr),
\end{equation}
where $\widehat{ES(\xi)_w^2}$ is defined as in (\ref{e:an-var-limit})
by replacing $f_{W_q}$ in $R_{q,w}$ and $f_W(w)$ by $\widehat f_{W_q}$,
and $z_\alpha$ is the $100\alpha$th percentile of $|{\mathcal
N}(0,1)|$. Another possibility is to use Proposition
\ref{p:an-estim-f-W-function} and set the $100\alpha\%$ confidence
``interval'' (set) across all $w$ simultaneously as
%
\begin{equation}\label{e:ci-2}
    \{f \dvtx \|\widehat f_W - f\|_{\infty,b} \leq N^{-1/2} \widehat
    q(\alpha)\}.
\end{equation}
Following \cite{grubelpitts1993} and
\cite{bogstedpitts2009}, the quantity $\widehat
q(\alpha)$ should be taken as the $100\alpha$th percentile of
%
\begin{equation}\label{e:R_N-hat}
    \widehat R_N(z) = \frac{1}{N^N} \sum_{(i_1,\ldots,i_N)\in
    \{1,\ldots,N\}^N} 1_{[0,z]}\Biggl(\sqrt{N} \Biggl\|
    S\Biggl(N^{-1}\sum_{k=1}^N 1_{\{W_{q,i_k} = \cdot\}}\Biggr) - S(\widehat f_{W_q})
    \Biggr\|_{\infty,b}\Biggr),
\end{equation}
which can be viewed as a bootstrapped version of the probability
%
\begin{equation}\label{e:R_N}
    R_N(z) = P\bigl(\sqrt{N} \|S(\widehat f_{W_q}) - S(f_{W_q})\|_{\infty,b}\leq z\bigr) =
P\bigl(\sqrt{N} \| \widehat f_W -  f_W \|_{\infty,b}\leq z\bigr).
\end{equation}
The choice of $\widehat q(\alpha)$ is justified if we can show that
$\widehat R_N$ and $R_N$ are asymptotically equivalent in distribution
in a suitable space of functions. This could be shown by following the
approach found in the proof of \cite{grubelpitts1993}, Proposition
3.15. The proof will not be given here. We should
also note that the discussion above assumes that the conditions of
Theorem \ref{t:an-estim-f-W} and Proposition
\ref{p:an-estim-f-W-function} are met, that is, $f_W(w)$ or
$f_{W_q}(w)$ has a sufficiently fast decay. In practice, these
conditions cannot be verified and caution should be exercised with the
resulting CI's. Another practical problem with Proposition
\ref{p:an-estim-f-W-function} and a related bootstrap procedure is that
$b$ is not known in advance. In the simulation study found in Section
\ref{s:simulations-number} below, we shall explore CI's given by
(\ref{e:ci-1}) and also use a~bootstrap procedure, but where
$\|\cdot\|_{\infty,b}$ is replaced by $\|\cdot\|_l$ with $\|x\|_l =
\sup_{1\leq n\leq l}|x_n|$ for some fixed $l$.

\subsection{Distribution of interrenewal times}
\label{s:times-sampled-to-renewals}

Here, we are interested in estimating the original distribution
function $F_D$ of interrenewal times. We focus on nonparametric
estimators, in the spirit of Buchmann and Gr\"{u}bel
\cite{buchmanngrubel2003} and B{\o}gsted and  Pitts
\cite{bogstedpitts2009}, and their basic properties (such as
asymptotic normality).

Estimation of $F_D$ is based on a theoretical inversion of the relation
(\ref{e:LT_sit}). (As mentioned in Section \ref{s:renewal-to-sampled},
forms of conditioning other than that on $W_q=s$ are possible and could
be dealt with by following the approach developed below.) Note that,
from (\ref{e:LT_sit}), we would expect
%
\begin{equation}\label{e:L_D-sampled->L_D}
    \widetilde F_D(v) = G_{a_s}(\widetilde F_{D_{q,i}|s}(v)),
\end{equation}
where $a_s = \{a_{s,n}\}_{n\in\Nset}$ is the reversion of the sequence
$A_s = \{A_{s,n}\}_{n\in\Nset}$ (see Section \ref{s:preliminaries}). In
terms of distribution functions, the relation
(\ref{e:L_D-sampled->L_D}) can be written as
%
\begin{equation}\label{e:F_D-sampled->F_D}
    F_D = \sum_{n=1}^\infty a_{s,n} F_{D_{q,i}|s}^{*n}.
\end{equation}

The following result provides sufficient conditions for
(\ref{e:F_D-sampled->F_D}) to hold and follows directly from~\cite{bogstedpitts2009}, Theorem 1. It uses
the following notation. Let $r(G_{A_s})$ and $r(G_{a_s})$ be the radii
of convergence of the corresponding power series $G_{A_s}$ and
$G_{a_s}$ (in the sense found in, e.g., \cite{henrici1974}).
By  B{\o}gsted and Pitts \cite{bogstedpitts2009}, Proposition 1,
there exist
\begin{equation}\label{e:sigma-A_s-a_s}
    \sigma(G_{A_s}) \in (0,r(G_{A_s})],\qquad
    \sigma(G_{a_s}) \in (0,r(G_{a_s})],
\end{equation}
such that $|w|<\sigma(G_{a_s})$ implies that there exists a unique $z$
where $|z|<\sigma(G_{A_s})$ and
\begin{equation}\label{e:A_s-a_s-identity}
        G_{A_s}(z) = w,\qquad z = G_{a_s}(w).
\end{equation}
%
%
We also need a suitable space of functions. For $\tau\geq 0$, let
$D_\tau[0,\infty) = \{f\dvtx[0,\infty)\mapsto \Rset \dvtx \|f\|_{\infty,\tau}:
= \sup_{t\geq 0} \mathrm{e}^{-\tau t} |f(t)|<\infty\}$.

\begin{theo}
[(\cite{bogstedpitts2009}, Theorem 1)]\label{t:f_D-smapled->F_D} Let $\tau
>0$ be such that
%
\begin{equation}\label{e:F_D-sampled->F_D-cond-1}
    \widetilde F_{D_{q,i}|s}(\tau) < \sigma(G_{a_s}).
\end{equation}
The series on the right-hand side of (\ref{e:F_D-sampled->F_D}) then
converges in $D_{\tau}[0,\infty)$. If, in addition,
%
\begin{equation}\label{e:F_D-sampled->F_D-cond-2}
    \widetilde F_D(\tau) < \sigma(G_{A_s}),
\end{equation}
then the relation (\ref{e:F_D-sampled->F_D}) holds.
\end{theo}

As discussed in \cite{bogstedpitts2009}, the
conditions (\ref{e:F_D-sampled->F_D-cond-1}) and
(\ref{e:F_D-sampled->F_D-cond-2}) always hold for large enough $\tau$.
We now turn to estimation based on the inversion
(\ref{e:F_D-sampled->F_D}). It is natural to consider the estimator
%
\begin{equation}\label{e:F_D-sampled->F_D-est}
    \widehat F_D = \sum_{n=1}^\infty \widehat a_{s,n} \widehat F_{D_{q,i}|s}^{*n}.
\end{equation}
Here, $\widehat a_s = \{\widehat a_{s,n}\}_{n\in\Nset}$ is the
reversion of the sequence $\widehat A_s = \{\widehat
A_{s,n}\}_{n\in\Nset}$ defined as
%
\begin{equation}\label{e:A_s-est}
    \widehat A_{s,n} = \frac{q^s}{\widehat f_{W_q} (s)} \sum_{w=s+n-1}^{\infty}
    \widehat f_W(w) {w-n\choose s-1}  (1-q)^{w-s},
\end{equation}
where $\widehat f_{W_q}$ and $\widehat f_W$ are given in
(\ref{e:estim-f_Wq}) and (\ref{e:estim-f_W-2}), respectively (i.e.,
$\widehat A_s$ is defined as in (\ref{eq:LT_sit_coef}), but where
$f_{W_q}$ and $f_W$ are replaced by their sample counterparts). For the
distribution function estimator, we take
\begin{equation}\label{e:F_D-sampled-est}
    \widehat F_{D_{q,i}|s}(t) = \frac{1}{N \widehat f_{W_q}(s)}
    \sum_{k=1}^N 1_{\{ D_{q,i,k} \leq t, W_{q,k}=s  \}},
\end{equation}
where $D_{q,i,k}$ are i.i.d. observations of $D_{q,i}$. An important
difference between (\ref{e:F_D-sampled->F_D-est}) and similar
estimators considered in
\cite{buchmanngrubel2003,bogstedpitts2009}
is that (\ref{e:F_D-sampled->F_D-est})
involves an estimator $\widehat a_s$ of $a_s$ (whereas $a_s$ was
assumed to be known in these related works). This obviously makes the
analysis of the estimator~(\ref{e:F_D-sampled->F_D-est}) more involved.
Also, note that the estimator $\widehat F_D$ in
(\ref{e:F_D-sampled->F_D-est}) is defined for fixed $s$ and $i$. In
addition, the estimator involves averages which are conditional on
$W_q=s$. In particular, the expressions (\ref{e:A_s-est}) and
(\ref{e:F_D-sampled-est}) are only defined for $\widehat f_{W_q}(s)\neq
0$.

We will show below that the estimator (\ref{e:F_D-sampled->F_D-est}) is
asymptotically normal under suitable assumptions. Note that the
estimator (\ref{e:F_D-sampled->F_D-est}) can be viewed as a functional
of $\widehat A_s$ (via $\widehat a_s$) and $\widehat F_{D_{q,i}|s}$. We
first need a suitable result on the asymptotic normality of the latter
quantities. Let $(\zeta = \{\zeta_n\}_{n\in\Nset},Z = \{Z(t)\}_{t\geq
0})$ be a zero-mean Gaussian process characterized by the following:
\begin{eqnarray} \label{e:an-A_s-F_D-sampled-zeta}
  \zeta_n &=& -\frac{\eta A_{s,n}}{f_{W_q}(s)} + \frac{q^s}{f_{W_q}(s)}
\sum_{w = n+s-1}^\infty S(\xi)_w {w-n\choose s-1} (1-q)^{w-s},\qquad
n\geq 1, \\\label{e:an-A_s-F_D-sampled-Z}
  Z(t) &=& -\frac{\eta F_{D_{q,i}|s}(t)}{f_{W_q}(s)} + \frac{1}{f_{W_q}(s)}
B(f_{W_q}(s) F_{D_{q,i}|s}(t)),\qquad t\geq 0,
\end{eqnarray}
where $S(\xi)$ appears in (\ref{e:an-estim-f-W}) of Theorem
\ref{t:an-estim-f-W}, $\eta$ is a zero-mean Gaussian variable with
$E\eta^2 = f_{W_q}(s) - f_{W_q}(s)^2$ and $B$ is a standard Brownian
bridge such that
\begin{eqnarray*}
 E\eta S(\xi)_w &=& 1_{\{s\geq w\}} {s\choose w}(-1)^{s-w} q^{-s} (1-q)^{s-w} f_{W_q}(s)
- f_{W_q}(s) f_{W}(w),\\
 E\eta B(f_{W_q}(s) F_{D_{q,i}|s}(t))  &=& f_{W_q}(s) \bigl(1 - f_{W_q}(s)\bigr) F_{D_{q,i}|s}(t), \\
  ES(\xi)_w B(f_{W_q}(s) F_{D_{q,i}|s}(t)) &=& 1_{\{s\geq w\}}
{s\choose w}\frac{(-1)^{s-w}}{q^{s}} (1-q)^{s-w} f_{W_q}(s) \bigl(1 -
f_{W_q}(s)\bigr)F_{D_{q,i}|s}(t).
\end{eqnarray*}
The proof of the following result can be found in Appendix
\ref{s:proofs}.

\begin{prop}\label{p:an-A_s-F_D-sampled}
Suppose that for some $z_0\geq 1$,
%
\begin{equation}\label{e:an-A_s-F_D-sampled-cond}
    \sum_{w=1}^\infty \sqrt{R_{q,w}} (1-q)^w z_0^w w^s < \infty,
\end{equation}
where $R_{q,w}$ is defined in (\ref{e:R_w-1}). Then,
\begin{equation}\label{e:an-A_s-F_D-sampled}
\sqrt{N} \bigl( (\widehat A_s,\widehat F_{D_{q,i}|s}) - (A_s,F_{D_{q,i}|s})
\bigr) \stackrel{d}{\rightarrow} (\zeta,Z)
\end{equation}
with the convergence in the space $(l_{\infty,z_0},D_0[0,\infty))$,
where the limit $(\zeta,Z)$ is characterized by
(\ref{e:an-A_s-F_D-sampled-zeta}) and (\ref{e:an-A_s-F_D-sampled-Z})
above.
\end{prop}

We are now ready to prove the asymptotic normality result for the
estimator $\widehat F_D$ in (\ref{e:F_D-sampled->F_D-est}). In addition
to the notation introduced before Theorem \ref{t:f_D-smapled->F_D}, we
shall also use the following. For a~formal power series $G(z) =
\sum_{n=0}^\infty a_n z^n$, let
%
\begin{equation}\label{e:nu_G}
\nu_G(\sigma) = \inf_{|z| = \sigma} |G(z)|.
\end{equation}
If $r(G)$ denotes the radius of convergence of $G$ as before, we also
let
\begin{equation}\label{e:r_0}
    r_0(G)= \sup\{\sigma_0\dvtx \sigma_0\leq r(G), \nu_G(\sigma)>0 \mbox{ for }
    0<\sigma<\sigma_0\}.
\end{equation}
The notation (\ref{e:nu_G})--(\ref{e:r_0}) follows that found in
\cite{henrici1974}; also recall the notation $\G_a(z)$ from
Section~\ref{s:preliminaries}.

\begin{theo}\label{t:an-F_D-sampled}
Suppose that the condition (\ref{e:an-A_s-F_D-sampled-cond}) holds with
$z_0\geq 1$. Also, suppose that for some $\tau>0$,
(\ref{e:F_D-sampled->F_D-cond-1}) and (\ref{e:F_D-sampled->F_D-cond-2})
hold, and
\begin{eqnarray}\label{e:an-F_D-sampled-cond-2}
{\G}_{a_s} (\widetilde F_{D_{q,i}|s}(\tau)) &\leq& z_0,\\\label{e:an-F_D-sampled-cond-3}
{\G}_{A_s}({\G}_{a_s}(\widetilde F_{D_{q,i}|s}(\tau))) &<& \max_{0<\sigma
<r_0(G_{A_s})\wedge z_0} \nu_{G_{A_s}}(\sigma).
\end{eqnarray}
Then,
\begin{equation}\label{e:an-F_D-sampled}
    \sqrt{N} (\widehat F_D - F_D) \stackrel{d}{\rightarrow} X
\end{equation}
with the convergence in the space $D_\tau[0,\infty)$. The limit $X$ is
a zero-mean Gaussian process which can be expressed as
\begin{equation}\label{e:an-F_D-sampled-limit}
    X = -\sum_{n=1}^\infty \bigl(a^{(1)}_s*(\zeta\circ a_s)\bigr)_n F_{D_{q,i}|s}^{*n}
    + \Biggl( \sum_{n=1}^\infty a_{s,n} n F_{D_{q,i}|s}^{*(n-1)} \Biggr) * Z,
\end{equation}
where $(\zeta,Z)$ is the limit process appearing in
(\ref{e:an-A_s-F_D-sampled}).
\end{theo}

\begin{rem}
The conditions
(\ref{e:an-F_D-sampled-cond-2}) and (\ref{e:an-F_D-sampled-cond-3}) hold
for large enough $\tau$. The presence of $z_0$ in~(\ref{e:an-A_s-F_D-sampled-cond}) and (\ref{e:an-F_D-sampled-cond-2})
is of technical interest: if $f_W(w)$ has faster decay, then $z_0$
could possibly be taken larger in (\ref{e:an-A_s-F_D-sampled-cond}) and
hence $\tau$ could be taken smaller in (\ref{e:an-F_D-sampled-cond-2}).
\end{rem}

\begin{pf*}{Proof of Theorem \ref{t:an-F_D-sampled}}
For notational simplicity, we will not write the index $s$ in
$A_s,a_s,\widehat A_s$ or $\widehat a_s$. By the Skorokhod
representation theorem and using Proposition
\ref{p:an-A_s-F_D-sampled}, we can suppose that
%
\begin{equation}\label{e:as-A_s-F_D-sampled}
\sqrt{N} \bigl( (\widehat A,\widehat F_{D_{q,i}|s}) - (A,F_{D_{q,i}|s}) \bigr)
\rightarrow (\zeta,Z)\qquad \mbox{a.s.}
\end{equation}
in the norm $(\|\cdot\|_{\infty,z_0},\|\cdot\|_{\infty,0})$. Write
\begin{eqnarray*}
 \sqrt{N} (\widehat F_D - F_D)  &=& \sqrt{N}
 \sum_{n=1}^\infty \widehat a_{n} \widehat F_{D_{q,i}|s}^{*n}
 - \sqrt{N} \sum_{n=1}^\infty a_{n} F_{D_{q,i}|s}^{*n} \\
   &=& \sqrt{N}
 \sum_{n=1}^\infty (\,\widehat a_{n} - a_{n}) (\widehat F_{D_{q,i}|s}^{*n}
 - F_{D_{q,i}|s}^{*n}) + \sqrt{N}
 \sum_{n=1}^\infty (\,\widehat a_{n} - a_{n}) F_{D_{q,i}|s}^{*n}  \\
   &&{} + \sqrt{N} \sum_{n=1}^\infty a_{n} (\widehat F_{D_{q,i}|s}^{*n}
 - F_{D_{q,i}|s}^{*n}) =: T_1 + T_2 + T_3.
\end{eqnarray*}
By B{\o}gsted and  Pitts \cite{bogstedpitts2009}, Theorem 2, we have
\[
T_3 \rightarrow \Biggl( \sum_{n=1}^\infty a_{n} nF_{D_{q,i}|s}^{*(n-1)} \Biggr)
* Z\qquad  \mbox{a.s.}
\]
in the norm $\|\cdot\|_{\infty,\tau}$, which is the second term in the
limit (\ref{e:an-F_D-sampled-limit}).

We next show that the term $T_2$ converges to the first term in the
limit (\ref{e:an-F_D-sampled-limit}). Observe that by using  \cite{buchmanngrubel2003}, Lemma 6(b) (see also the
inequality used in the proof of their Lemma 7), we
have\vspace*{-1pt}
\begin{eqnarray*}
 \|E_2\|_{\infty,\tau}  &: =&  \Biggl\| \sqrt{N} \sum_{n=1}^\infty (\widehat a_{n} - a_{n})
F_{D_{q,i}|s}^{*n} - \sum_{n=1}^\infty \bigl(-a^{(1)}*(\zeta\circ
a)\bigr)_n F_{D_{q,i}|s}^{*n} \Biggr\|_{\infty,\tau} \\
   &\leq&  \sum_{n=1}^\infty \bigl|\sqrt{N}  (\,\widehat a_{n} - a_{n}) +
\bigl(a^{(1)}*(\zeta\circ a)\bigr)_n \bigr| \|  F_{D_{q,i}|s}^{*n}
\|_{\infty,\tau}\\
  &\leq& C \sum_{n=1}^\infty \bigl|\sqrt{N}  (\,\widehat a_{n} - a_{n}) +
\bigl(a^{(1)}*(\zeta\circ a)\bigr)_n \bigr| \widetilde F_{D_{q,i}|s}(\tau)^n.
\end{eqnarray*}
Writing $\sqrt{N}  (\,\widehat a - a) + a^{(1)}*(\zeta \circ a) =
(\sqrt{N}  (\,\widehat a\circ A - I) + (a^{(1)}\circ A)
* \zeta ) \circ a$ and using the inequality
$\G_{x\circ y}(z) \leq \G_x(\G_y(z))$, we further obtain that\vspace*{-1pt}
\begin{eqnarray*}
  \|E_2\|_{\infty,\tau} &\leq& C \sum_{n=1}^\infty \bigl| \bigl(\sqrt{N} (\,\widehat
a\circ A - I) + \bigl(a^{(1)}\circ A\bigr)
* \zeta \bigr)_n\bigr| (\G_{a}(\widetilde F_{D_{q,i}|s}(\tau)))^n \\
   &\leq &C(R_1 + R_2 + R_3 + R_4),
\end{eqnarray*}
where, with $z_1 = \G_{a}(\widetilde F_{D_{q,i}|s}(\tau))$,
\begin{eqnarray*}
  R_1 &=&  \sum_{n=1}^\infty \bigl| \bigl(\sqrt{N} \bigl(\widehat a\circ (\widehat A + A - \widehat A) -
  \widehat a\circ \widehat A\,\bigr) - \bigl(\widehat a^{(1)}\circ \widehat A\,\bigr)
* \bigl(\sqrt{N}(A - \widehat A) \bigr) \bigr)_n\bigr| z_1^n, \\
  R_2 &=& \sum_{n=1}^\infty \bigl| \bigl( \bigl(\widehat a^{\,(1)}\circ \widehat A\,\bigr)
* \bigl(\sqrt{N}(A - \widehat A) \bigr) - \bigl( a^{(1)}\circ \widehat A\,\bigr)
* \bigl(\sqrt{N}(A - \widehat A) \bigr) \bigr)_n\bigr| z_1^n,   \\
  R_3 &=& \sum_{n=1}^\infty \bigl| \bigl( \bigl( a^{(1)}\circ \widehat A\,\bigr)
* \bigl(\sqrt{N}(A - \widehat A) \bigr) - \bigl( a^{(1)}\circ A\bigr)
* \bigl(\sqrt{N}(A - \widehat A) \bigr) \bigr)_n\bigr| z_1^n,  \\
  R_4 &=& \sum_{n=1}^\infty \bigl| \bigl(\bigl ( a^{(1)}\circ A\bigr)
* \bigl(\sqrt{N}(\widehat A - A) \bigr) - \bigl(a^{(1)}\circ A\bigr)
* \zeta \bigr)_n\bigr| z_1^n.
\end{eqnarray*}
We will next show that $R_k\to 0$ a.s., $k=1,2,3,4$.

For the term $R_1$, by using Proposition \ref{p:first-order} in
Appendix \ref{s:appendix-bounds-taylor}, we first observe that\vspace*{-2pt}
\begin{eqnarray*}
  R_1 &\leq & \frac{\sqrt{N}}{2} \sum_{n=1}^\infty \bigl(
    \bigl(\widehat a^{\,(2)} \circ\bigl(\widehat A_+ + (A-\widehat A)_+ \bigr) \bigr)
    * (A-\widehat A)_+*(A-\widehat A)_+
  \bigr)_n z_1^n  \\
   &=& \frac{1}{2\sqrt{N}} \G_{\widehat a^{\,(2)}}\bigl(\G_{\widehat A}(z_1) +
   \G_{A-\widehat A}(z_1)\bigr) \bigl(\G_{\sqrt{N} (A-\widehat A)}(z_1)\bigr)^2.
\end{eqnarray*}
By (\ref{e:as-A_s-F_D-sampled}), and since $z_1 \leq z_0$ by the
assumption (\ref{e:an-F_D-sampled-cond-2}), we have $\G_{\sqrt{N}
(A-\widehat A)}(z_1) \to \G_{\zeta}(z_1)$ a.s. Next, we want to show
that\vspace*{-2pt}
\begin{equation}\label{e:G-a^2-bound}
\G_{\widehat a^{\,(2)}}\bigl(\G_{\widehat A}(z_1) +
   \G_{A-\widehat A}(z_1)\bigr) \leq C\qquad \mbox{a.s.}
\end{equation}
for some random constant $C$. For this, we further examine the radius
of convergence of $\G_{\widehat a^{\,(2)}}$.  By the Cauchy--Hadamard
formula (see, e.g., \cite{henrici1974}, Theorem 2.2a, page
77),
\[
r\bigl(\G_{\widehat a^{\,(2)}}\bigr) = r(\G_{\widehat a}) = r(G_{\widehat a}).
\]
By applying the inequality after the proof of \cite{henrici1974}, Theorem 2.2b, page
99, we then have
\[
r\bigl(\G_{\widehat a^{\,(2)}}\bigr)  \geq \max_{0<\sigma <r_0(G_{\widehat
A})}\nu_{G_{\widehat A}}(\sigma) \geq \max_{0<\sigma <r_0(G_{\widehat
A}) \wedge z_0}\nu_{G_{\widehat A}}(\sigma),
\]
where the notation $r_0(G)$, $\nu_G$ was introduced above. Since
$G_{\widehat A}(z)$ converges to $G_{A}(z)$ a.s. and uniformly on
$|z|\leq z_0$, we have
\[
\max_{0<\sigma <r_0(G_{\widehat A}) \wedge z_0}\nu_{G_{\widehat
A}}(\sigma) \to \max_{0<\sigma <r_0(G_{A}) \wedge
z_0}\nu_{G_{A}}(\sigma)\qquad \mbox{a.s.}
\]
The relation (\ref{e:G-a^2-bound}) now follows from the relations
above, the fact that $\G_{\widehat A}(z_1) + \G_{A-\widehat A}(z_1) \to
\G_A(z_1)$ a.s. and the assumption (\ref{e:an-F_D-sampled-cond-3}). We
can now conclude that $R_1\to 0$ a.s.

For the term $R_2$, note that for fixed $K\geq 1$,
\begin{eqnarray*}
  R_2 &\leq & \sum_{n=1}^\infty \bigl| \bigl( \widehat a^{\,(1)}  - a^{(1)} \bigr)_n\bigr| (\G_{\widehat
  A}(z_1))^n \G_{\sqrt{N}(A - \widehat A)}(z_1) \\
   &\leq & \Biggl(\sum_{n=1}^K \bigl| \bigl( \widehat a^{\,(1)}  - a^{(1)} \bigr)_n\bigr| (\G_{\widehat
  A}(z_1))^n + \sum_{n=K+1}^\infty \bigl|\widehat a^{\,(1)}_n\bigr| (\G_{\widehat
  A}(z_1))^n   \\
   &&\phantom{\Biggl(}{} +  \sum_{n=K+1}^\infty \bigl|a^{(1)}_n\bigr| (\G_{\widehat
  A}(z_1))^n \Biggr)  \G_{\sqrt{N}(A - \widehat A)}(z_1)  =: R_{2,1} + R_{2,2} + R_{2,3}.
\end{eqnarray*}
For a fixed $K$, $R_{2,1}\to 0$ a.s. by using the a.s. convergence of
$\widehat A$ in (\ref{e:as-A_s-F_D-sampled}). On the other hand,
arguing as for the term $R_1$ above, we can make sure that $R_{2,2}$
and $R_{2,3}$ are arbitrarily small for large enough $K$. For the term
$R_3$, by using Proposition \ref{p:first-order},
\begin{eqnarray*}
  R_3 &\leq & \sum_{n=1}^\infty \bigl| \bigl( a^{(1)}\circ \widehat A
- a^{(1)}\circ A\bigr)_n \bigr| z_1^n \G_{\sqrt{N}(A - \widehat A)}(z_1)
 \\
   &\leq &  \sum_{n=1}^\infty \bigl(
    a_+^{(2)}\circ\bigl(A_+ + (\widehat A-A)_+\bigr) *(\widehat A-A)_+
   \bigr)_n  z_1^n \G_{\sqrt{N}(A - \widehat A)}(z_1)\\
   &= & \frac{1}{\sqrt{N}} \G_{a^{(2)}}\bigl(\G_{A}(z_1) + \G_{\widehat
   A-A}(z_1)\bigr) \bigl(\G_{\sqrt{N}(A - \widehat A)}(z_1)\bigr)^2
\end{eqnarray*}
and we have $R_3\to 0$ by arguing as for the term $R_1$. By using the
bound
\[
R_4\leq \G_{a^{(1)}}(\G_A(z_1)) \G_{\sqrt{N}(\widehat A - A) - \zeta
}(z_1),
\]
we get that $R_4\to 0$ a.s.

Finally, we shall prove that the term $T_1$ is asymptotically
negligible. Note that, as in the proof of
\cite{buchmanngrubel2003}, Proposition 8,
\begin{eqnarray*}
  \| T_1 \|_{\infty,\tau} &\leq & \sqrt{N}
 \sum_{n=1}^\infty |\widehat a_{n} - a_{n}| \| \widehat F_{D_{q,i}|s}^{*n}
 - F_{D_{q,i}|s}^{*n}\|_{\infty,\tau} \\
   &\leq & \bigl\|\sqrt{N}( \widehat F_{D_{q,i}|s} - F_{D_{q,i}|s}) \bigr\|_{\infty,0} \sum_{n=1}^\infty |\widehat a_{n} -
a_{n}|
 \widetilde H_n(\tau),
\end{eqnarray*}
where
\[
H_n = \sum_{j=0}^{n-1} \widehat F_{D_{q,i}|s}^{\,*j} *
F_{D_{q,i}|s}^{*(n-1-j)}.
\]
Using the fact that $\,\widetilde{\!{\widehat{F}}\hspace*{0.5pt}}_{D_{q,i}|s}(\tau) \to
\widetilde  F_{D_{q,i}|s}(\tau)$ a.s., we can further deduce that
\[
\| T_1 \|_{\infty,\tau} \leq C\bigl\|\sqrt{N}( \widehat F_{D_{q,i}|s}
 - F_{D_{q,i}|s}) \bigr\|_{\infty,0} \sum_{n=1}^\infty |\widehat a_{n} - a_{n}|
 \bigl((1+\delta)\widetilde F_{D_{q,i}|s}(\tau)\bigr)^n
\]
for any arbitrarily small $\delta>0$ and  all $N$ large enough. The
convergence of the latter series to $0$ can be proven by arguing as for
the term $T_2$ and using the fact that $\delta>0$ can be chosen
arbitrarily small.
\end{pf*}

Despite its obvious theoretical interest, Theorem
\ref{t:an-F_D-sampled} is of limited relevance in practice where the
conditions of the theorem cannot be verified and $\tau$ is unknown.
Note also that the form of the limiting Gaussian process in
(\ref{e:an-F_D-sampled-limit}) is not easily tractable. In a simulation
study found in Section~\ref{s:simulations-times} below,\vadjust{\goodbreak} we  make
qualitative observations on the performance of $\widehat F_D$ and the
use of bootstrap for CI's. The CI's provided below will be based on the
sup-norm on an interval $[0,T]$, with fixed $T$. Justification of the
bootstrap procedure is beyond the scope of this paper.

\section{Simulation study}
\label{s:simulation}

Here, we provide a simulation study for the inference of the
distributions of the number of renewals and interrenewal times.

\subsection{Simulations for number of renewals}
\label{s:simulations-number}

As mentioned in \hyperref[s:intro]{Introduction}, the estimator $\widehat f_W$ in
(\ref{e:estim-f_W-1}) is known to perform poorly as $q$ gets smaller.
Here, we shall re-examine this fact and the performance of the
estimator $\widehat f_W$ via the asymptotic variance $ES(\xi)^2_w$ in
(\ref{e:an-var-limit}) of the normal limit in (\ref{e:an-estim-f-W}).
In this regard, two regimes should be distinguished:\vspace*{-1pt}
\begin{eqnarray}\label{e:two-regimes}
     & \bullet &\mbox{ stable }\Bigl(\sup_{w\in\Nset} R_{q,w} < \infty \mbox{ or } \sup_{w\in\Nset} ES(\xi)_w^2 <
     \infty\Bigr);\nonumber\\ [-8pt]\\ [-8pt]
    &\bullet& \mbox{ explosive }  \Bigl(\sup_{w\in\Nset} R_{q,w} = \infty \mbox{ or }
 \sup_{w\in\Nset}ES(\xi)_w^2 = \infty\Bigr),\nonumber
\end{eqnarray}
where $R_{q,w}$ is defined in (\ref{e:R_w-1}), that is,
\[
R_{q,w} = \sum_{s=w}^\infty {s\choose w}^2
\frac{(1-q)^{2(s-w)}}{q^{2s}} f_{W_q}(s).
\]
As will be seen below, the performance of the estimator $\widehat f_W$
largely depends on which of the two regimes is considered. Before
providing some simulations, we take a closer look at these regimes for
the distributions $f_W$ considered in Section \ref{s:examples}.

\subsubsection*{Geometric distribution} Suppose that $f_W$ follows a
geometric distribution with parameter $c\in (0,1),$ as in Section
\ref{s:examples}. We shall derive lower and upper bounds for $R_{q,w}$.
For the lower bound, observe that
\[
f_{W_q}(s) \geq \sum_{w=s}^\infty q^s (1-q)^{w-s} f_W(w) = C_1 (cq)^s
\]
for some constant $C_1$ (which depends on $c$ and $q$) and hence that
\begin{equation}\label{e:R-lower-bound-geom}
    R_{q,w} \geq C_1 \sum_{s=w}^\infty \frac{(1-q)^{2(s-w)}}{q^{2s}}
    ( qc)^s = C_2 \biggl( \frac{c}{q} \biggr)^w
\end{equation}
for some constant $C_2$ (which depends on $c$ and $q$). For the upper
bound of $R_{q,w}$, we need a bound on binomial coefficients. We shall
use the standard (but rough) bound given by ${n\choose k} \leq n^k/k!$.
We shall also use the following auxiliary result, proved in Appendix
\ref{s:proofs}.

\begin{lem}\label{l:bound-integral}
Let $a\in (0,1)$. Then, for some constant $C>0$ and all $s\geq 1$,
\begin{eqnarray}\label{e:bound-integral}
  \int_s^\infty w^s a^w\,\mathrm{d}w  &\leq & C s^s a^s \bigl( 1 + (\ln a^{-1})^{-s}\bigr), \\\label{e:bound-integral2}
  \int_s^\infty w^{2s} a^w\,\mathrm{d}w &\leq & C s^{2s} a^s \biggl( 1 + \biggl(\frac{1}{2}\ln
  a^{-1}\biggr)^{-2s}\biggr).
\end{eqnarray}
\end{lem}


Using the lemma, with generic constants $C_k$,
\begin{eqnarray*}
  f_{W_q}(s) & \leq & \frac{C_1}{s!} q^s (1-q)^{-s} \sum_{w=s}^\infty w^s \bigl(c(1-q)\bigr)^w \\
   &\leq &  \frac{C_2}{s!} q^s (1-q)^{-s} \int_s^\infty w^s \bigl(c(1-q)\bigr)^w\,\mathrm{d}w \\
   &\leq &  \frac{C_3}{s!} q^s (1-q)^{-s} s^s \bigl(c(1-q)\bigr)^s \bigl(1 + \bigl(\ln \bigl(c(1-q)\bigr)^{-1}\bigr)^{-s}\bigr) \\
   &\leq &  C_4 (\mathrm{e} c q)^s \bigl(1 + \bigl(\ln \bigl(c(1-q)\bigr)^{-1}\bigr)^{-s}\bigr),
\end{eqnarray*}
where the last inequality follows from  Stirling's formula. Hence, by
arguing similarly and using~(\ref{e:bound-integral2}), we
have
\begin{eqnarray}  \label{e:R-upper-bound-geom}
   R_{q,w} &\leq & C_1 \frac{(1-q)^{-2w}}{(w!)^2} \sum_{s=w}^\infty s^{2w} \frac{(1-q)^{2s}}{q^{2s}}
   \biggl( (\mathrm{e}cq)^s + \biggl( \frac{\mathrm{e}cq}{\ln (c(1-q))^{-1}} \biggr)^s \biggr) \nonumber \\
   &\leq & C_2 \frac{(1-q)^{-2w}}{(w!)^2} w^{2w} \biggl(
    \frac{(1-q)^{2w}}{q^{2w}} (\mathrm{e}cq)^w \biggl( 1 + \biggl(\frac{1}{2}\ln \frac{q}{(1-q)^2 \mathrm{e}c} \biggr)^{-2w} \biggr)\nonumber \\
   & &\hphantom{C_2 \frac{(1-q)^{-2w}}{(w!)^2} w^{2w} \biggl(} {}+ \frac{(1-q)^{2w}}{q^{2w}} \biggl( \frac{\mathrm{e}cq}{\ln (c(1-q))^{-1}}
   \biggr)^w\nonumber\\
    &&\hphantom{C_2 \frac{(1-q)^{-2w}}{(w!)^2} w^{2w} \biggl({}+\hspace*{1.5pt}}{}\times \biggl( 1
    + \biggl( \frac{1}{2}\ln \frac{q\ln (c(1-q))^{-1}}{(1-q)^2 \mathrm{e}c} \biggr)^{-2w} \biggr)\biggr)
     \nonumber\\ [-8pt]\\ [-8pt]
   & \leq & C_3  \biggl( \frac{\mathrm{e}^3c}{q} \biggr)^w  \biggl(
    1 + \biggl(\ln \frac{1}{c(1-q)}\biggr)^{-w} + \biggl(\frac{1}{2} \ln \frac{q}{(1-q)^2 \mathrm{e}c}
    \biggr)^{-2w}\nonumber \\
   & &\hphantom{C_3  \biggl( \frac{\mathrm{e}^3c}{q} \biggr)^w  \biggl(}{}+ \biggl( \frac{1}{2}\ln \frac{q\ln (c(1-q))^{-1}}{(1-q)^2 \mathrm{e}c}
   \biggr)^{-2w}\nonumber\\
   &&\hphantom{C_3  \biggl( \frac{\mathrm{e}^3c}{q} \biggr)^w  \biggl(}{} + \biggl(\ln \frac{1}{c(1-q)}\biggr)^{-w}
   \biggl( \frac{1}{2}\ln \frac{q\ln (c(1-q))^{-1}}{(1-q)^2 \mathrm{e}c} \biggr)^{-2w} \biggr) \nonumber\\
   &\leq & C_4 \biggl( \frac{\mathrm{e}^3c}{q} \biggr)^w,\qquad \mbox{if }
   \frac{\mathrm{e}^3c}{q} \leq 1,\nonumber
\end{eqnarray}
where, for the last inequality, we used the fact that all  three log
terms (before raising to the powers $-w$ and $-2w$) are bigger than 1
under $\mathrm{e}^3c/q \leq 1$.

The bounds (\ref{e:R-lower-bound-geom}) and
(\ref{e:R-upper-bound-geom}) show that, for the geometric distribution
$f_W$,
\begin{eqnarray}\label{e:regimes-geom}
    \frac{c}{q} \leq \mathrm{e}^{-3} \quad &\Rightarrow& \quad \mbox{stable
    regime},\nonumber\\ [-8pt]\\ [-8pt]
    \frac{c}{q} > 1 \quad &\Rightarrow& \quad \mbox{explosive regime}.\nonumber
\end{eqnarray}
What happens in the range $c/q \in (\mathrm{e}^{-3},1]$ remains an open
question. Our experience in practice suggests that the critical point
for $c/q$ is closer to $1$ and that $\mathrm{e}^{-3}$ is a very rough bound.

\subsubsection*{Heavy-tailed distribution}
Suppose that $f_W$ follows a
heavy-tailed distribution with parameter $\alpha$ (Section
\ref{s:examples}). We will show that in this case, the distribution is
always in the explosive regime (and this happens for all $\alpha>0$,
not just $\alpha\in(1,2)$). Indeed, observe that
\[
f_{W_q}(s) \geq  C \sum_{w=s}^\infty q^s (1-q)^{w-s} w^{-\alpha-1} \geq
C q^s s^{-\alpha-1}
\]
and hence
\begin{equation}\label{e:R-lower-bound-heavy}
 R_{q,w} \geq  C \sum_{s=w}^\infty \frac{(1-q)^{2(s-w)}}{q^{2s}} q^s s^{-\alpha-1}
   \geq  C \frac{w^{-\alpha-1}}{q^w}.
\end{equation}
It remains to observe that the lower bound in
(\ref{e:R-lower-bound-heavy}) diverges as $w\to \infty$ since
$q\in(0,1)$.

\begin{rem}\label{r:explosive}
In the examples above, note that in the explosive regime, $R_{q,w}$ or
$ES(\xi)_w^2$ diverges at an exponential rate. The term ``explosive''
was chosen to reflect this fact. Also, note  that in the explosive
regime, in order to have convergence of $\widehat f_W(w)$ for all $w$,
we naturally need to consider weighted spaces, as in Proposition
\ref{p:an-estim-f-W-function} above.
\end{rem}

\begin{rem}[(The impact of small $\bolds{q}$)]
Note that $R_{q,w}\uparrow \infty$ as $q\downarrow 0$. Since
$R_{q,w}\geq q^{-2w}f_{W_q}(w)$, when $w$ is larger and $q$ is small,
$R_{q,w}$ (or $\operatorname{Var}(\widehat f_W(w))$) can also be too large for
practical purposes (e.g., with $w=10$ and $q=0.1$, $R_{q,w}\geq 10^{20}
f_{W_q}(w)$). Moreover, as with the geometric distribution above, we
may expect that most of distributions of interest (with unbounded
support) belong to the explosive regime for small enough $q$. These
observations reiterate the current understanding that inference of
$f_{W_q}$ becomes impractical for small $q$.
\end{rem}

\begin{figure}

\includegraphics{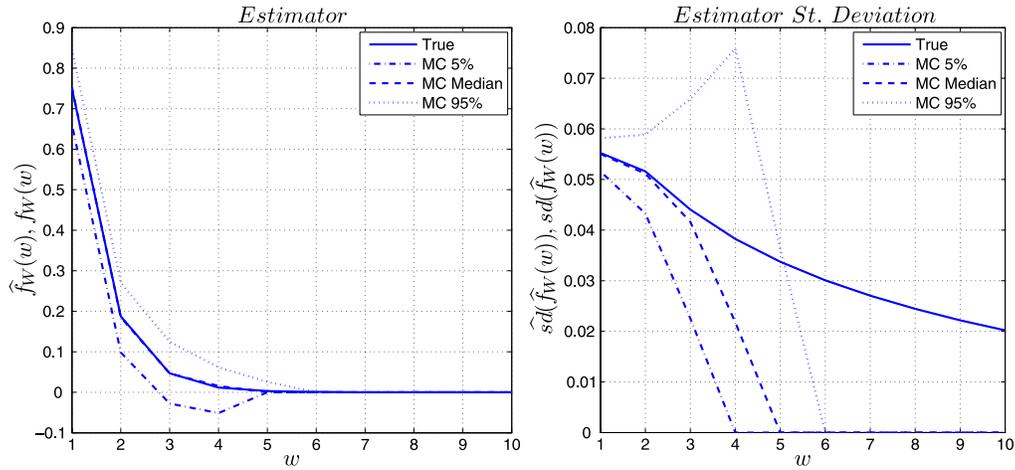}

\caption{Performance of the estimator $\widehat
f_W(w)$ for geometric distribution with $c=0.25$, $q=0.6$, $N=500$.
\textup{Left plot:} true $f_W(w)$ and MC-based $5\%$, $50\%$, $95\%$ percentiles
of $\widehat f_W(w)$. \textup{Right plot:} true asymptotic standard deviation
$\mathit{sd}(\widehat f_W(w))$ and MC-based $5\%$, $50\%$, $95\%$ percentiles of
$\widehat{\mathit{sd}}(\widehat f_W(w))$.}\label{fig:geometric}
\end{figure}

\begin{figure}[b]

\includegraphics{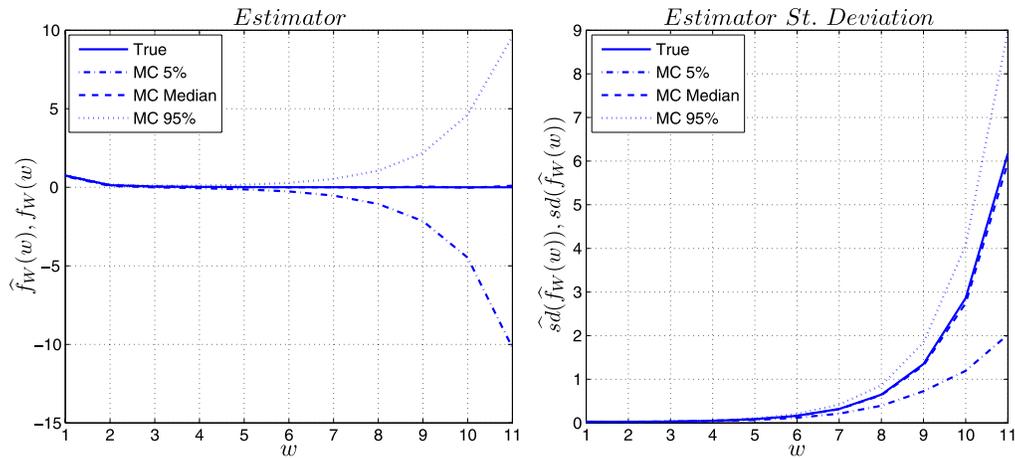}

\caption{Performance of the estimator $\widehat
f_W(w)$ for Pareto distribution with $\alpha=1.5$, $q=0.7$, $N=1000$.
\textup{Left plot:} true $f_W(w)$ and MC-based $5\%$, $50\%$, $95\%$ percentiles
of $\widehat f_W(w)$. \textup{Right plot:} true asymptotic standard deviation
$\mathit{sd}(\widehat f_W(w))$ and MC-based $5\%$, $50\%$, $95\%$ percentiles of
$\widehat{\mathit{sd}}(\widehat f_W(w))$.}\label{fig:pareto}
\end{figure}

In Figures \ref{fig:geometric} and \ref{fig:pareto}, we illustrate the
performance of the estimator $\widehat f_W(w)$ in the two regimes
above. Figure \ref{fig:geometric} is for the geometric distribution
with $c=0.25$, $q=0.6$, $N=500$, and Figure~\ref{fig:pareto} is for the
Pareto distribution with $\alpha = 1.5$, $q=0.7$, $N=1000$. The
left-hand plots in these figures depict the true p.m.f. $f_W(w)$, and
the $5\%$, $50\%$ and $95\%$ percentiles of the distribution of the
estimator $\widehat f_W(w)$, based on $1000$ Monte Carlo (MC)
replications. The right-hand plots in these figures depict the true
asymptotic standard deviation $(ES(\xi)^2_w/N)^{1/2}$ or $\mathit{sd}(\widehat
f_W(w))$, and the $5\%$, $50\%$ and $95\%$ percentiles of the
distribution of the estimator $(\widehat{ES(\xi)^2_w}/N)^{1/2}$ or
$\widehat{\mathit{sd}}(\widehat f_W(w))$, again based on 1000 MC replications.

Note that by the definition (\ref{e:estim-f_W-2}) of $\widehat f_W(w)$,
we have $\widehat f_W(w)=0$ for all $w$ larger than the size of the
largest sampled flow(s). The largest sampled flow(s) in the study
corresponding to Figure~\ref{fig:geometric} consists of $5$ sampled
packets. For this reason, the $\widehat f_W(w)$'s are all zero for
$w\geq 6$ in Figure \ref{fig:geometric} (left plot) and obviously have
zero estimated standard deviation for $w\geq 6$ in Figure
\ref{fig:geometric} (right plot). Figure  \ref{fig:pareto}, on the
other hand, does not depict zero $\widehat f_W(w)$'s because the size
of the largest sampled flow(s) is much greater. Indeed, only a short
range of $w=1,\ldots,11$ is considered in Figure \ref{fig:pareto}.  As
noted in Remark \ref{r:explosive} above, the estimator standard
deviation grows exponentially with increasing~$w$. For larger $w$, the
features of both plots in Figure \ref{fig:pareto} would be dominated by
this exponential blow-up, even if there is a drop to zero in $\widehat
f_W(w)$ and its estimated standard deviation after the largest sampled
flow(s).

Figure  \ref{fig:pareto} depicts typical estimation in the explosive
regime: $\widehat f_W(w)$ and its estimated standard deviation blow up
(before dropping to zero). In Figure  \ref{fig:geometric}, on the other
hand, $\widehat f_W(w)$ and its estimated standard deviation are stable
and reach zero. This is the situation associated with the stable
regime. Thus, even though the parameters $c$ and $q$ in Figure
\ref{fig:geometric} do not satisfy the sufficient condition for the
stable regime in (\ref{e:regimes-geom}), we still expect that these
parameters actually lead to the stable regime. (Recall that, as
indicated earlier, the condition in (\ref{e:regimes-geom}) is likely to
be very strong.) Also, observe that Figure \ref{fig:pareto} is for a
sample size of only $N=1000$, and roughly the first $5$ frequencies can
be estimated with confidence. With larger sample sizes, for example,
with $N=10\thinspace000$ and $N=100\thinspace000$, roughly the first
$10$ and $12$ frequencies, respectively, can be estimated with
confidence.

\begin{figure}

\includegraphics{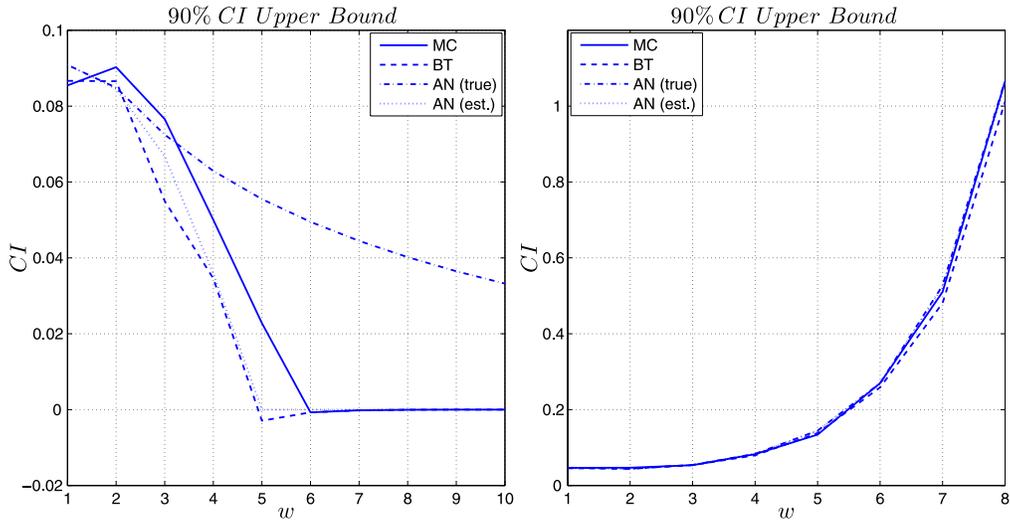}

\caption{$95\%$ CI upper bounds for $\widehat f_W(w) -
f_W(w)$, that is, for $\widehat f_W(w)$ centered at true $f_W(w)$, for
the two situations of Figures \protect\ref{fig:geometric} and \protect\ref{fig:pareto}.
The left plot corresponds to Figure \protect\ref{fig:geometric} and the right
plot corresponds to Figure \protect\ref{fig:pareto}. Both plots consist of the
$95\%$ MC-based percentile of $\widehat f_W(w)-f_W(w)$ (MC), median
bound using bootstrap (BT) and bounds based on the asymptotic normality
(AN) result (\protect\ref{e:an-estim-f-W}) for $\widehat f_W(w)$ with true
(``true'') and estimated (``est.'') asymptotic variance.}\label{fig:ci}
\end{figure}

In Figure \ref{fig:ci}, we give an idea of the appropriateness of
confidence intervals (CI's). The left- and right-hand plots of Figure
\ref{fig:ci} correspond to the situations considered in Figures
\ref{fig:geometric} and \ref{fig:pareto}, respectively. The $95\%$
MC-based percentile of $\widehat f_W(w)$ is depicted as in Figures
\ref{fig:geometric} and \ref{fig:pareto}. This percentile is, in fact,
centered at the true $f_W(w)$. It can be thought of as an ideal upper
$90\%$ CI bound, corresponding to dealing with an exact distribution of
$\widehat f_W(w)$. We also depict analogous bounds based on
bootstrapping (BT) and the asymptotic normality (AN) result
(\ref{e:an-estim-f-W}) for $\widehat f_W(w)$, with either true
(``true'') or estimated (``est.'') asymptotic variance. In the cases of
BT and AN (est.) bounds, what are depicted are, in fact, the median
bounds obtained from MC realizations. As with MC, the BT, AN (true) and
AN (est.) bounds are centered at $f_W(w)$. Note from Figure~\ref{fig:ci} that CI's based on both BT and AN (est.) approaches, which
are the only practical alternatives, are quite satisfactory.

\subsection{Simulations for interrenewal times}
\label{s:simulations-times}

We present here a simulation study for the estimator $\widehat F_D$ in
(\ref{e:F_D-sampled->F_D-est}) of the distribution of interrenewal
times. Two cases are considered. In case 1, we take geometric $f_W$
with $c=0.25$ and $q=0.6$, as in Section \ref{s:simulations-number}
(corresponding to the stable regime). In case 2, geometric $f_W$ is
taken with $c=0.7$ and $q=0.6$ (corresponding to the explosive regime;
the simulation results are similar if the Pareto distribution is used
as in Section \ref{s:simulations-number}). In both cases, $D$ is
supposed to be an exponential random variable with parameter $1$. The
respective sample sizes are $N=500$ and $N=1000$.

\begin{figure}

\includegraphics{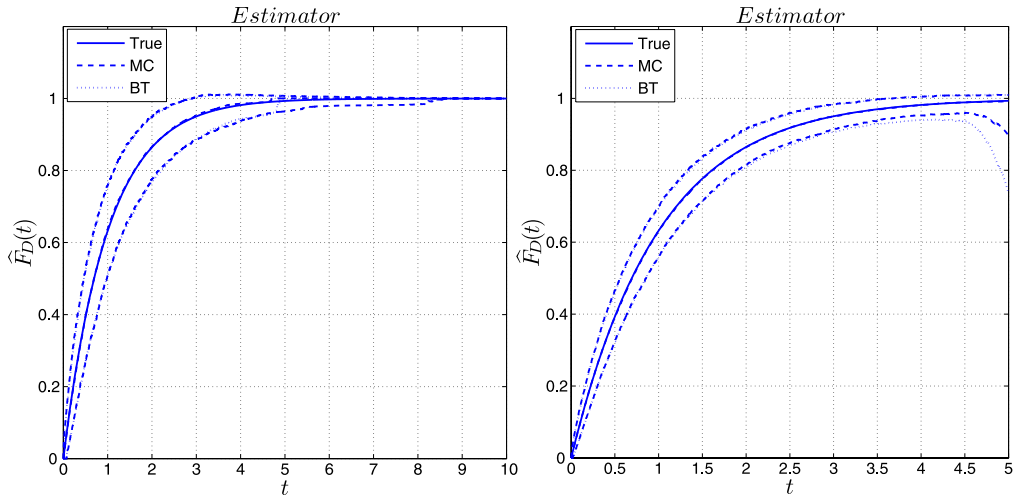}

\caption{Performance of $\widehat F_D$ and the
bootstrap CI's for the exponential distribution with parameter $1$.
\textup{Left plot:} $N=500$, $q=0.6$ and geometric $f_W$ with $c=0.25$. \textup{Right
plot:} $N=500$, $q=0.6$ and geometric $f_W$ with $c=0.7$. Plots include
true $F_D$, MC-based $5\%$, $50\%$ and $95\%$ percentiles of $\widehat
F_D$, and the medians of $90\%$ bootstrap (BT) CI's.}\label{fig:times}
\end{figure}

Figure \ref{fig:times} presents simulation results for case 1 in the
left-hand plot and case 2 in the right-hand plot. We depict the true
function $F_D$ and the $5\%$, $50\%$ and $95\%$ percentiles of the
distribution of $\widehat F_D$, based on 1000 MC realizations. The
plots also contain the medians of $90\%$ bootstrap-based CI's, computed
from 1000 MC replications. The latter suggests, in particular, that
such CI's are quite appropriate. We should also note that the estimator
$\widehat F_D$ is here based only on the first sampled interrenewal
times when $s=2$, that is, $i=1$ and $s=2$ in
(\ref{e:F_D-sampled->F_D-est}). In addition, the estimator is computed
by truncating the infinite sum in its definition
(\ref{e:F_D-sampled->F_D-est}) and by evaluating convolutions
numerically through discretizing convolution integrals.

Note, from Figure \ref{fig:times}, that estimation is satisfactory,
even for case 2 which corresponds to the explosive regime. This seems
quite surprising and needs further explanation. In fact, in the setting
of case 2, only the first several values of $\widehat f_W(w)$ are
really relevant for the estimator $\widehat F_D$. This can be seen from
the definition (\ref{e:A_s-est}) of the sequence $\widehat A_s$, which
enters into $\widehat F_D$ after reversion. Note, from
(\ref{e:A_s-est}), that in calculating $\widehat A_{s,n}$, the first
series term $\widehat f_W(n+1)$ is weighted by $(1-q)^{n-1}$. This
additional factor thus helps to keep the blow-up of $\widehat f_W$
under control. Moreover, from the examples above, if the growth of
$\widehat f_W(w)$ is thought to be of  order $q^{-w}$, then the factor
$(1-q)^w$ annihilates this growth when $q>0.5$. In practice, and in
case 2, we also observe that, even though $\widehat f_W(w)$ is highly
varying for larger $w$, the sequence $\widehat A_{s,n}$ decays to zero
rapidly. Since the first few values of $\widehat f_W(w)$ can be
estimated with confidence even in the explosive regime, the above
explains the satisfactory performance of $\widehat F_D$ in case 2 of
Figure \ref{fig:times}. On the other hand, note also that the $t$-range
is smaller in the right-hand plot of Figure \ref{fig:times} and that
the variance of the estimator starts to diverge at the boundary $T=5$
of the range. In fact, this divergence would be more pronounced (and
would dominate the plot) if we increased the range of $t$. Thus, the
performance of the estimator is satisfactory, but only to some time
point $T$.

\begin{figure}[b]

\includegraphics{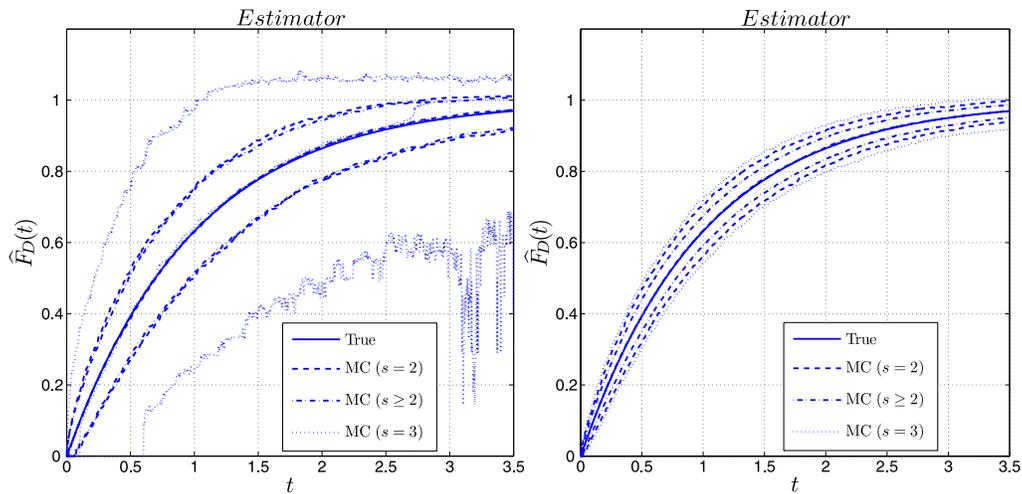}

\caption{Performance of the estimator $\widehat
F_D(t)$ for two cases from Figure \protect\ref{fig:times}. Plots include
MC-based $5\%$, $50\%$ and $95\%$ percentiles of $\widehat F_D(t)$ when
conditioning on $W_q=s$ with $s=2$, $s\geq 2$ and $s=3$.}\label{fig:s=2-3->2}
\end{figure}

In Figure \ref{fig:s=2-3->2}, we also illustrate the performance of
$\widehat F_D$ when using several other forms of conditioning on $W_q$.
The two plots in the figure correspond to the respective plots of
Figure \ref{fig:times}. They depict the $5\%$, $50\%$ and $95\%$
percentiles of $\widehat F_D$ when conditioning on $W_q=s$ with $s=2$,
$s\geq 2$ and $s=3$. (In the case of conditioning on $s\geq 2$, the
formulas (\ref{LT:geq-s}) and (\ref{series-coef:geq-s}) are used.) The
$5$--$95\%$ interpercentile range is the smallest when conditioning on
$s\geq 2$, although comparable to that when $s=2$ is used. An obvious
note to add here is that different types of conditioning (for the same
$N$) lead to different effective sample sizes of the data used in the
computation of $\widehat F_D$, with that for $s\geq 2$ being the
largest and that for $s=3$ being the smallest. The interpercentile
range is much larger for $s=3$ than for $s=2$ in the left-hand plot
because the difference between the corresponding effective sample sizes
is much larger for case 1.

\begin{appendix}\label{appendix}

\section{Proofs}
\label{s:proofs}

\begin{pf*}{Proof of Theorem \ref{t:LT_joint_sit}} For $1\leq i_1<
\cdots <i_n< s$ and $t_1\geq 0,\ldots,t_n \geq 0$, we have
\begin{eqnarray}\label{p:LT_joint_sit_aux1}
&&\hspace*{-20pt}P(D_{q,i_1}\leq t_1,\ldots,D_{q,i_n} \leq t_n,W_{q}=s) \nonumber \\ [-8pt]\\ [-8pt]
&&\hspace*{-20pt}\quad=  \sum_{w=s}^{\infty} f_W(w) P(W_{q}=s|W=w)
P(D_{q,i_1}\leq t_1,\ldots,D_{q,i_n} \leq t_n|W=w,W_{q}=s).\nonumber
\end{eqnarray}
Let $M_j$, $j=1,\ldots,n$, be the number of renewals not sampled
between the $i_{j}$th and $(i_{j}+1)$th sampled renewals plus 1. If
$W=w$ and $W_q=s$ are fixed and $M_1=m_1,\ldots,M_{j-1}=m_{j-1}$, then
$M_j$ can take the values $1,\ldots,w-m_1-\cdots - m_{j-1} -(s-j)$ and
hence
\begin{eqnarray}\label{p:LT_joint_sit_aux2}
&&P(D_{q,i_1} \leq t_1,\ldots,D_{q,i_n} \leq t_n|W=w,W_{q}=s)\nonumber\\
&&\quad=\sum_{m_1=1}^{w-(s-1)}\sum_{m_2=1}^{w-m_1-(s-2)}\ldots
\sum_{m_n=1}^{w-m_1-\cdots -
m_{n-1} -(s-n)}
F_D^{*m_1}(t_1) \cdots F_D^{*m_n}(t_n)\\
&&\qquad{}\times P(M_{1}=m_1,\ldots,M_{n}=m_n|W=w,W_{q}=s).\nonumber
\end{eqnarray}

We shall next derive an expression for the probability
$P(M_{1}=m_1,\ldots,M_{n}=m_n|W=w,W_{q}=s)$ in
(\ref{p:LT_joint_sit_aux2}). The indices $i_1,\ldots,i_n$ are not
necessarily consecutive, but still form separate blocks of consecutive
indices (e.g., $i_1=1$, $i_2=2$, $i_3=4$, $i_4=8$, $i_5=9$ form three
separate blocks of indices $\{i_1,i_2\}$, $\{i_3\}$ and $\{i_4,i_5\}$).
We need additional notation to keep track of   these separate blocks
and the  gaps between them. Thus, let $1=j_1<j_2<\cdots< j_k \leq n$
denote subindices $j$ in $i_j$ for which the corresponding $i_j$ is the
start of a separate block (e.g., with $i_1=1$, $i_2=2$, $i_3=4$,
$i_4=8$, $i_5=9$ above, we have $j_1=1$, $j_2=3$, $j_3=4$). With this
notation, note that $r_u=j_{u+1}-j_u$, $1\leq u<k-1$ and
$r_k=n-{j_k}+1$ denote the sizes of the $k$ separate blocks. Also, let
$x_0=i_{j_1}-1$, $x_u=i_{j_{u+1}}- i_{j_{u}} - r_{u}-1$, $1\leq u<k-1$,
and $x_k=s-i_{j_k}-r_k$ be the numbers of sampled renewals,
respectively, before the $i_{j_1}\!$th sampled renewal, between the
$(i_{j_{u}} +r_{u})$th and $i_{j_{u+1}}\!$th sampled renewals, and after
the $(i_{j_k}+r_{k})$th sampled renewal.

Given $W=w$ and $W_q =s$, the number of possible distinct locations of
the sampled renewals satisfying $M_1=m_1, \ldots,M_n=m_n$ can be
obtained by considering the  possible initial locations of the sampled
renewals $i_{j_u}$, along with the number of possible distinct
locations of the sampled renewals before  the $i_{j_1}$th, between  the
$(i_{j_{u}} +r_{u})$th and $i_{j_{u+1}}$th, and after the
$(i_{j_k}+r_{k})$th sampled renewal, that is,
\begin{eqnarray*}
&&\sum_{l_1=x_0+1}^{w-\widehat m_1 -\cdots -\widehat m_k-x_1-\cdots
-x_k-(k-1)} {l_1-1\choose x_0} \sum_{l_2=l_1+\widehat m_1+x_1
+1}^{w-\widehat m_2 -\cdots -\widehat m_k-x_2-\cdots -x_k-(k-1)}
{l_2-l_1-\widehat m_1-1\choose x_1}\\
&&\quad{}\times \cdots
\times \sum_{l_{{k-1}}=l_{k-2}+\widehat m_{k-2}+x_{k-1}+1}^{w-\widehat
m_{k-1} -\widehat m_k-x_{k-1} -x_k-1}
{l_{k-1}-l_{k-2}-\widehat m_{k-2}-1\choose x_{k-2}}\\
&&\quad{}\times \sum_{l_k=l_{k-1}+\widehat m_{k-1}+x_{k-1}+1}^{w-\widehat
m_k-x_k} {l_k - l_{k-1} -\widehat{m}_{k-1}-1\choose x_{k-1}} {w-
l_k - \widehat m_k\choose x_{k}},
\end{eqnarray*}
where $l_u$ is the location of the $i_{j_u}$th sampled\vspace*{1pt} renewal and
$\widehat m_u =\sum_{v=j_u}^{j_u + r_u-1} m_{v}$ is the number of
renewals between the $i_{j_u}$th and $(i_{j_{u}} +r_{u})$th sampled
renewals  plus 1. By making the change of variables
$l_u^{\prime}=l_u-\sum_{v=1}^{u-1} \widehat m_v$, $u=2,\ldots,k$, the
last expression becomes
\begin{eqnarray*}
&&\sum_{l_1=x_0+1}^{w-m-x_1-\cdots-x_{k}-(k-1)} {l_1-1\choose x_0}
\sum_{l_2^{\prime}=l_1^{\prime}+x_1+1}^{w-m-x_2-\cdots-x_{k}-(k-2)}
{l_2^{\prime}-l_1^{\prime}-1\choose x_1}\\
&&\quad{} \times \cdots \times
\sum_{l_{k-1}^{\prime}=l_{k-2}+x_{k-1}+1}^{w-m-x_{k-1}-x_{k}-1}
{l_{k-1}^{\prime}-l_{k-2}^{\prime}-1\choose x_{k-2}}\\
 &&\quad{}\times\sum_{l_k^{\prime}=l_{k-1}^{\prime}+x_{k-1}+1}^{w-m-x_{k}} {l_k^{\prime} -
l_{k-1}^{\prime}-1\choose x_{k-1}} {w-m-
l_k^{\prime}\choose x_{k}}={w-m\choose x_0+\cdots+x_{k}+k},
\end{eqnarray*}
where the last identity can be viewed as a generalization of the
identity (\ref{eq:ind}) used in the proof of Theorem \ref{t:LT_sit}.
Since $x_0+\cdots+x_{k}+k=s-n$, we deduce that
\begin{equation} \label{p:LT_joint_sit_aux3}
P(M_{1}=m_1,\ldots,M_{n}=m_n|W=w,W_{q}=s)
 ={w- m\choose s-n}\bigg/{w\choose s}.
\end{equation}

By using (\ref{p:LT_joint_sit_aux1})--(\ref{p:LT_joint_sit_aux3}), we
have that
\begin{eqnarray}\label{e:DF_joint_sit}
F_{\mathbf D_{q,n}|s}(t)
&=& \frac{P(D_{q,i_1}\leq t_1,\ldots,D_{q,i_n} \leq t_n,W_{q}=s)}{f_{W_q} (s)}\nonumber \\
&=&\frac{1}{f_{W_q} (s)} \sum_{w=s}^{\infty} f_W(w) q^{s}
(1-q)^{w-s}\nonumber\\
&&\hphantom{\frac{1}{f_{W_q} (s)} \sum_{w=s}^{\infty}}{}\times\sum_{m_1=1}^{w-(s-1)} \cdots \sum_{m_n=1}^{w-m_1-\cdots-m_{n-1}-(s-n)}
{w-m\choose s-n}\\
&&\hphantom{\frac{1}{f_{W_q} (s)} \sum_{w=s}^{\infty}}{}\times F_D^{*m_1} (t_1) \cdots F_D^{*m_n} (t_n) \nonumber \\
&=&\sum_{m_1=1}^{\infty} \ldots \sum_{m_n=1}^{\infty}
B_{s,(m_1,\ldots,m_n)} F_D^{*m_1} (t_1) \cdots F_D^{*m_n} (t_n),\nonumber
\end{eqnarray}
where $B_{s,\mathbf m}$ is given by (\ref{e:LT_joint_sit_coef}). The
relation (\ref{e:LT_joint_sit-1}) is obtained by taking the LT in
(\ref{e:DF_joint_sit}).
\end{pf*}

\begin{pf*}{Proof of Proposition \ref{p:LT_s_V}} For $s\geq 2$ and
$t\geq 0$, we have
%
\begin{equation} \label{p:LT_s_V_aux1}
 P(V_q \leq t, W_q=s) =   \sum_{w=2}^{\infty} f_W(w) P(W_q=s|W=w) P(V_q \leq t|W=w,W_q=s).
\end{equation}
Now, letting $M$ be the number of original renewals between the first
and  last sampled renewals plus 1, we have
\begin{eqnarray}\label{p:LT_s_V_aux2}
 P(V_q \leq t|W=w,W_q=s) &=& \sum_{m=1}^{w-1} F_D^{*m} (t) P(M=m|W=w,W_s=s) \nonumber
 \\ [-8pt]\\ [-8pt]
&=& \sum_{m=1}^{w-1} F_D^{*m}(t) (w-m){m-1\choose s-2}\bigg/{w\choose s},\nonumber
\end{eqnarray}
where $(w-m)$ is the number of possible distinct locations of the first
sampled renewal and ${m-1\choose s-2}$ is the number of possible
locations of sampled renewals between the first and last sampled
renewals when $M=m$. By using (\ref{p:LT_s_V_aux1}) and
(\ref{p:LT_s_V_aux2}), we deduce that
\begin{equation} \label{e:DF_s_V}
F_{V_q|2^+}(t)  = \frac{1}{P(W_q \geq 2)} \sum_{s=2}^{\infty} P(V_q
\leq t, W_q=s)
= \sum_{m=1}^{\infty} C_{m}  F_D^{*m} (t),
\end{equation}
which can also be written as in (\ref{e:LT_s_V}).
\end{pf*}

\begin{pf*}{Proof of Proposition \ref{p:f_Wq->f_W-1}} The equality
in (\ref{e:f_Wq->f_W-1-cond-main}) follows from (\ref{e:pmf-s_W}) and
\begin{eqnarray*}
&&\sum_{s=n}^\infty {s\choose n} \frac{(1-q)^{s-n}}{q^s}
\sum_{w=s}^\infty {w\choose s} q^{s} (1-q)^{w-s} f_{W}(w)\\
&&\quad =
\sum_{w=n}^\infty f_{W}(w) (1-q)^{w-n} \sum_{s=n}^w {s\choose n}
{w\choose s}\\
&&\quad= \sum_{w=n}^\infty f_{W}(w) (1-q)^{w-n} {w\choose n} \sum_{s=n}^w
{w-n\choose s-n} = \sum_{w=n}^\infty f_{W}(w) (1-q)^{w-n} {w\choose n}
2^{w-n}.
\end{eqnarray*}
Then, by the assumption (\ref{e:f_Wq->f_W-1-cond-main}), we can
substitute (\ref{e:pmf-s_W}) into the right-hand side of
(\ref{e:f_Wq->f_W}) and use  Fubini's theorem to obtain
\begin{eqnarray*}
   & &\sum_{s=w}^\infty {s\choose w} \frac{(-1)^{s-w}}{q^s} (1-q)^{s-w}
\sum_{n=s}^\infty {n\choose s} q^{s} (1-q)^{n-s} f_{W}(n) \\
  &&\quad= \sum_{n=w}^\infty f_W(n) (1-q)^{n-w} \sum_{s=w}^n
                        {s\choose w} {n\choose s} (-1)^{s-w}   \\
   &&\quad=  \sum_{n=w}^\infty f_W(n) (1-q)^{n-w}  {n\choose w}  \sum_{s=w}^n
                       {n-w\choose s-w} (-1)^{s-w} \\
   &&\quad= \sum_{n=w}^\infty f_W(n) (1-q)^{n-w}  {n\choose w}
   \sum_{k=0}^{n-w} {n-w\choose k} (-1)^{k} = f_W(w),
\end{eqnarray*}
where we have used the identity $\sum_{k=0}^K {K\choose k} (-1)^k = 0$
if $K\geq 1$, and $=1$ if $K=0$.
\end{pf*}

 \begin{pf*}{Proof of Proposition \ref{p:T-S}} We consider only the
case $l=2$ (and $z_2=0$). Note that
\begin{eqnarray*}
  T^{(2)}(x)_n &=&  \sum_{i=n}^\infty {i\choose n} T^{(1)}(x)_i (- z_1)^{i-n}\\
   &=&  \sum_{i=n}^\infty {i\choose n} (- z_1)^{i-n} \sum_{k=i}^\infty {k\choose i} \frac{x_k}{q^k} \bigl(z_1 -
    (1-q)\bigr)^{k-i} \\
  &=& \sum_{k=n}^\infty \frac{x_k}{q^k}\sum_{i=n}^k {i\choose n}{k\choose i} (- z_1)^{i-n} \bigl(z_1 -
    (1-q)\bigr)^{k-i}  \\
  &=& \sum_{k=n}^\infty \frac{x_k}{q^k} {k\choose n}\sum_{i=n}^k {k-n\choose i-n} (- z_1)^{i-n} \bigl(z_1 -
    (1-q)\bigr)^{k-i}  \\
  &=& \sum_{k=n}^\infty \frac{x_k}{q^k} {k\choose n} \bigl(-(1-q)\bigr)^{k-n}
  = S(x)_n.
\end{eqnarray*}
The change of the order of summation above can be justified by using
Fubini's theorem and  the assumption (\ref{e:T-S-cond}).\
\end{pf*}

\begin{pf*}{Proof of Proposition \ref{p:f_Wq->f_W-2}} The proof is
similar to that of Propositions \ref{p:f_Wq->f_W-1} and \ref{p:T-S},
and only the case $l=2$ (and $z_2=0$) will be considered. By using
(\ref{e:pmf-s_W}), we observe that
%
\begin{eqnarray}\label{e:T^1}
  T^{(1)}(f_{W_q})_n &=&  \sum_{i=n}^\infty {i\choose n}  \frac{(z_1 - (1-q))^{i-n}}{q^i}
  \sum_{k=i}^\infty {k\choose i} q^i (1-q)^{k-i} f_W(k)  \nonumber \\
   &=&  \sum_{k=n}^\infty f_W(k) \sum_{i=n}^k {i\choose n} {k\choose i}
     \bigl(z_1 - (1-q)\bigr)^{i-n} (1-q)^{k-i}\nonumber \\ [-8pt]\\ [-8pt]
   &=& \sum_{k=n}^\infty f_W(k){k\choose n} \sum_{i=n}^k  {k-n\choose i-n}
     \bigl(z_1 - (1-q)\bigr)^{i-n} (1-q)^{k-i}  \nonumber \\
   &=& \sum_{k=n}^\infty f_W(k) {k\choose n} z_1^{k-n},\nonumber
\end{eqnarray}
where the order of the summation above is changed using Fubini's
theorem. Indeed, the application of  Fubini's theorem is possible as
long as $|z_1 - (1-q)| + (1-q)<1$ or $z_1\in C_0$, which is one of the
conditions on $z_1$. By using (\ref{e:T^1}), we further get that
\begin{eqnarray*}
  T^{(2)}(f_{W_q})_n &=& \sum_{i=n}^\infty {i\choose n} (-z_1)^{i-n} \sum_{k=i}^\infty
  {k\choose i} z_1^{k-i} f_W(k) \\
   &=& \sum_{k=n}^\infty f_W(k) \sum_{i=n}^k {i\choose n}
  {k\choose i} (-z_1)^{i-n} z_1^{k-i} = f_W(n),
\end{eqnarray*}
where the use of  Fubini's theorem is justified by the fact that
$|{-}z_1|+|z_1|<1$ or $z_2 = 0\in C_1$, which is one of the conditions on
$z_2$.
\end{pf*}

\begin{pf*}{Proof of Proposition \ref{p:necerssary}} Write
$T^{(1)}(\xi^{(1)}) = T^{(1)}(\xi^{(1)}_1) + T^{(1)}(\xi^{(1)}_2)$,
where
\[
\xi^{(1)}_{1,s} = \cases{
\xi^{(1)}_s, &\quad if $s\leq j$,\cr
0, &\quad if $s> j$,
}\qquad
\xi^{(1)}_{2,s} = \cases{ 0, &\quad
if $s\leq  j$, \cr
\xi^{(1)}_s, &\quad if $s> j$
}
\]
and observe that, for fixed $j$, $T^{(1)}(\xi^{(1)}_1)$ and
$T^{(1)}(\xi^{(1)}_2)$ are independent. We can then also write
\begin{equation}\label{e:T^2_xi1_j-decomp}
T^{(2)}\bigl(\xi^{(1)}\bigr)_{j,w} = T^{(2)}\bigl(\xi^{(1)}_1\bigr)_{j,w} +
T^{(2)}\bigl(\xi^{(1)}_2\bigr)_{j,w} = T^{(2)}\bigl(\xi^{(1)}_1\bigr)_{w} +
T^{(2)}\bigl(\xi^{(1)}_2\bigr)_{j,w},
\end{equation}
where $T^{(2)}(\xi^{(1)}_1)_{w}$ and $T^{(2)}(\xi^{(1)}_2)_{j,w}$ are
independent for fixed $j$.  By Proposition \ref{p:T-S}, we have that
\[
T^{(2)}\bigl(\xi^{(1)}_1\bigr)_{w} = S\bigl(\xi^{(1)}_1\bigr)_w = \sum_{s=w}^j {s\choose w}
\frac{(-1)^{s-w}}{q^s} (1-q)^{s-w} f_{W_q}(s)^{1/2} \eta_s.
\]
By using the fact that the two terms in (\ref{e:T^2_xi1_j-decomp}) are
Gaussian and independent, we can then write
\begin{equation}\label{e:T^2_xi1_j-decomp-2}
E\mathrm{e}^{\mathrm{i}\theta T^{(2)}(\xi^{(1)})_{j,w}} = \mathrm{e}^{-(1/2)\theta^2
\sigma^2_{1,w}} \mathrm{e}^{-(1/2)\theta^2 \sigma^2_{2,w}},
\end{equation}
where
\begin{eqnarray*}
\sigma^2_{1,w} &=& E\bigl(T^{(2)}\bigl(\xi^{(1)}_1\bigr)_{w} \bigr)^2 =
   \sum_{s=w}^j {{s\choose w}}^2 \frac{(1-q)^{2(s-w)}}{q^{2s}}
   f_{W_q}(s),\\[-2pt]
\sigma^2_{2,w}&=& E\bigl(T^{(2)}\bigl(\xi^{(1)}_2\bigr)_{j,w} \bigr)^2.
\end{eqnarray*}
In view of (\ref{e:T^2_xi1_j-decomp-2}), for $T^{(2)}(\xi^{(1)})_{j,w}$
to converge in distribution, it is necessary that the series
$\sigma^2_{1,w}$ converges, that is, that
(\ref{e:an-estim-f-W-cond-main}) holds. (In other words, if
(\ref{e:an-estim-f-W-cond-main}) does not hold, then the limit of~(\ref{e:T^2_xi1_j-decomp-2}) is zero.)
\end{pf*}

\begin{pf*}{Proof of Proposition \ref{p:an-estim-f-W-function}} The
proof is in two steps: (a) $\sqrt{N}(\widehat f_{W_q} - f_{W_q}) \to
_d \xi$ in the space $l_{\infty,q^{-1}(b'+1-q)}$, where $\xi$ is the
process appearing in the proof\vspace*{1pt} of Theorem \ref{t:an-estim-f-W}; (b)
the mapping $S$ in (\ref{e:S_x}) is continuous from
$l_{\infty,q^{-1}(b'+1-q)}$ to $l_{\infty,b}$.

For (a), it is enough to show the tightness of $\sqrt{N}(\widehat
f_{W_q} - f_{W_q})$ in $l_{\infty,q^{-1}(b'+1-q)}$. By the results of
\cite{vakhaniatarieladzechobayan1987}, page
229 (see also
\cite{buchmann2001}, Theorem 2.5), this follows from the fact that
for any $\epsilon>0$,
\begin{equation}\label{e:an-estim-f-W-function-tightness}
    \lim_{m\to \infty} \sup_{N\geq 1}
    P\Bigl(\sup_{s>m} \bigl(q^{-1}(b'+1-q)\bigr)^s  \bigl|\sqrt{N}\bigl(\widehat f_{W_q}(s) - f_{W_q}(s)\bigr) \bigr|
    > \epsilon\Bigr) = 0.
\end{equation}
For the latter, observe that
\begin{eqnarray*}
&&\epsilon P\biggl(\sup_{s>m} \biggl ( \frac{b'+1-q}{q} \biggr)^s \bigl|\sqrt{N}\bigl(\widehat
f_{W_q}(s) - f_{W_q}(s)\bigr) \bigr|
    > \epsilon\biggr)\\[-2pt]
&&\quad\leq \sum_{s=m+1}^\infty  \biggl( \frac{b'+1-q}{q} \biggr)^s  E \bigl|\sqrt{N}\bigl(\widehat
f_{W_q}(s) - f_{W_q}(s)\bigr) \bigr|\\[-2pt]
&&\quad\leq \sum_{s=m+1}^\infty  \biggl( \frac{b'+1-q}{q} \biggr)^s  \sqrt{E
\bigl|\sqrt{N}\bigl(\widehat f_{W_q}(s) - f_{W_q}(s)\bigr) \bigr|^2} \leq
\sum_{s=m+1}^\infty  \biggl( \frac{b'+1-q}{q} \biggr)^s f_{W_q}(s)^{1/2}
\end{eqnarray*}
and that (\ref{e:an-estim-f-W-function-tightness}) follows from the
assumption (\ref{e:an-estim-f-W-function-cond}).

For (b), observe that
\begin{eqnarray*}
\|S(x) - S(y)\|_{\infty,b} &=& \sup_{n\geq 0} b^n |S(x)_n - S(y)_n| \leq
\sum_{n=0}^\infty b^n |S(x)_n - S(y)_n|\\[-2pt]
&\leq& \sum_{n=0}^\infty b^n \sum_{i=n}^\infty {i\choose n}
\frac{(1-q)^{i-n}}{q^i} |x_i-y_i| = \sum_{i=0}^\infty |x_i - y_i| \biggl(
\frac{b+1-q}{q} \biggr)^i\\[-2pt]
&\leq& \sup_{n\geq 0} \biggl( \frac{b' +1-q}{q} \biggr)^n |x_n - y_n|
\sum_{i=0}^\infty \biggl( \frac{b+1-q}{b'+1-q} \biggr)^i\\[-2pt]
& =&
\|x-y\|_{\infty,q^{-1}(b'+1-q)} \sum_{i=0}^\infty \biggl(
\frac{b+1-q}{b'+1-q} \biggr)^i.
\end{eqnarray*}
\vspace*{-12pt}
\upqed\end{pf*}

\begin{pf*}{Proof of Proposition \ref{p:an-A_s-F_D-sampled}} With
$X^N = \sqrt{N} (\widehat A_s - A_s)$ and $Y^N = \sqrt{N} (\widehat
F_{D_{q,i}|s} - F_{D_{q,i}|s})$, it is enough to show that (a) the
finite-dimensional distributions of $(X^N,Y^N)$ converge to those of
$(\zeta,Z)$; (b) the sequence $(X^N,Y^N)$ is tight.
The latter condition follows from the following two conditions: (b1)
the sequence $X^N = \{X^N_n\}_{n\in\Nset}$ is tight; (b2) the
sequence $Y^N$ is tight.
The tightness of $Y_N$ in (b2) is a standard result,  part of the
empirical central limit theorem (see, e.g., \cite{pollard1984}, Section
V.2) and will not be proven here. By the results of
\cite{vakhaniatarieladzechobayan1987}, the condition (b1) can be
replaced by (b1) for any $\epsilon>0$, $\lim_{m\to \infty}
\sup_{N\geq 1} P(\sup_{n>m} z_0^n |X^N_n| > \epsilon) = 0$.

For the latter (b1), it is enough to show the same relation, but
where $X^N$ is replaced by $\widetilde X^N$, which is defined without
$q^s/\widehat f_{W_q}(s)$ and $q^s/f_{W_q}(s)$ in $\widehat A_s$
and $A_s$, respectively. With the sequence $\widetilde X^N$, observe
that
\begin{eqnarray*}
   \epsilon P\Bigl(\sup_{n>m} z_0^n |\widetilde X^N_n| > \epsilon\Bigr) &\leq & \sum_{n=m+1}^\infty  z_0^n E|\widetilde X^N_n| \\
   &\leq & \sum_{n=m+1}^\infty \sum_{w=n+s-1}^\infty  z_0^n E\bigl|\sqrt{N} \bigl(\widehat f_W(w) - f_W(w)\bigr)\bigr|
   {w-n\choose s-1} (1-q)^{w-s} \\
   &=& \sum_{w=m+s}^\infty E\bigl|\sqrt{N} \bigl(\widehat f_W(w) - f_W(w)\bigr)\bigr| (1-q)^{w-s} \sum_{n=m+1}^{w-s+1} {w-n\choose s-1}z_0^k  \\
   &\leq & \sum_{w=m+s}^\infty E\bigl|\sqrt{N} \bigl(\widehat f_W(w) - f_W(w)\bigr)\bigr|
   (1-q)^{w-s} z_0^{w-s+1} {w-m\choose s},
\end{eqnarray*}
where we have used the fact that $z_0\geq 1$ and $\sum_{n=m+1}^{w-s+1}
{w-n\choose s-1} = {w-m\choose s}$. Since $E|\sqrt{N} (\widehat f_W(w)
- f_W(w))|^2\leq R_{q,w}$, as in the proof of Theorem
\ref{t:an-estim-f-W}, we  further get that
\begin{eqnarray*}
   \epsilon P\Bigl(\sup_{n>m} z_0^n |\widetilde X^N_n| > \epsilon\Bigr)
   & \leq & \sum_{w=m+s}^\infty \sqrt{ E\bigl|\sqrt{N} \bigl(\widehat f_W(w) - f_W(w)\bigr)\bigr|^2}
   (1-q)^{w-s} z_0^{w-s+1}
   {w-m\choose s} \\
   & \leq & \sum_{w=m+s}^\infty \sqrt{R_{q,w}} (1-q)^{w-s} z_0^{w-s+1}
   {w-m\choose s}.
\end{eqnarray*}
The condition (b1) now follows from the assumption
(\ref{e:an-A_s-F_D-sampled-cond}).

For the convergence of the finite-dimensional distributions in (a),
we show only the convergence of $X^N_n$ (the general case can be
considered along similar lines). For $j\geq n+s-1$, define
\[
X^N_{j,n} = \sqrt{N} (\widehat A^j_{s,n} - A^j_{s,n}),
\]
where
\[
 \widehat A^j_{s,n} = \frac{q^s}{\widehat f_{W_q} (s)} \sum_{w=n+s-1}^{j}
    \widehat f_W(w) {w-n\choose s-1}  (1-q)^{w-s}
\]
and similarly with $A^j_{s,n}$, and also
\[
\zeta^j_n = -\frac{\eta A_{s,n}}{f_{W_q}(s)} + \frac{q^s}{f_{W_q}(s)}
\sum_{w = n+s-1}^j S(\xi)_w {w-n\choose s-1} (1-q)^{w-s}.
\]
It is enough to show that:
\begin{longlist}
    \item $X^N_{j,n} \stackrel{d}{\to} \zeta^j_n$ as $N\to\infty$;
    \item $\zeta^j_n \stackrel{d}{\to} \zeta_n$ as $j\to\infty$;
    \item for any $\delta>0$, $\limsup_{j\to \infty} \limsup_{N\to \infty} P(|X^N_{j,n} -
X^N_n|> \delta ) = 0$.
\end{longlist}
The convergence in (i) follows directly from Theorem
\ref{t:an-estim-f-W}. The convergence in (ii) follows using
$E|S(\xi)_w| \leq \sqrt{ES(\xi)_w^2} \leq \sqrt{R_{q,w}}$ and the
assumption (\ref{e:an-A_s-F_D-sampled-cond}). Condition (iii) can be
proven as in part (b1) above using the assumption
(\ref{e:an-A_s-F_D-sampled-cond}).
\end{pf*}

\begin{pf*}{Proof of Lemma \ref{l:bound-integral}} Using
integration by parts, we obtain
\begin{eqnarray*}
   \int_s^\infty w^s a^w\,\mathrm{d}w & = &  \frac{s^s a^s}{\ln a^{-1}} +
 \frac{s}{\ln a^{-1}} \int_s^\infty w^{s-1} a^w\,\mathrm{d}w\\
   &\leq & \frac{s^s a^s}{\ln a^{-1}} +\frac{s^s a^s}{(\ln a^{-1})^2} +
 \frac{s^2}{(\ln a^{-1})^2} \int_s^\infty w^{s-2} a^w\,\mathrm{d}w \\
   &\leq & s^s a^s \sum_{k=1}^s \frac{1}{(\ln a^{-1})^k} +
   \frac{s^s}{(\ln a^{-1})^s} \int_s^\infty a^w\,\mathrm{d}w \\
   &=& s^s a^s \sum_{k=1}^{s+1} \frac{1}{(\ln a^{-1})^k}
   =  \frac{s^s a^s}{\ln a^{-1}}\cdot \frac{(\ln a^{-1})^{-s-1} - 1}{(\ln a^{-1})^{-1} - 1},
\end{eqnarray*}
from which the bound (\ref{e:bound-integral}) follows. The inequality
(\ref{e:bound-integral2}) can be proven similarly.\
\end{pf*}

\section{Bounds on remainder terms in Taylor expansions of compositions of
power series} \label{s:appendix-bounds-taylor}

Here, we prove a result which was used several times in the proof of
Theorem \ref{t:an-F_D-sampled}.

\begin{prop}\label{p:first-order}
For any formal power series $x(z)=\sum_{n=1}^{\infty}x_n z^n$,
$y(z)=\sum_{n=1}^{\infty}y_n z^n$, $\epsilon (z)=\sum_{n=1}^{\infty}
\epsilon_n z^n$ and $n\geq 1$, we have
\begin{eqnarray}\label{e:first-order}
    \bigl|\bigl(x\circ (y+\epsilon)  - x\circ y\bigr)_n\bigr| &\leq&  \bigl(\bigl(x^{(1)}_+ \circ (y_+ + \epsilon_+)\bigr)*
\epsilon_+  \bigr)_n,
\\\label{e:second-order}
    \bigl|\bigl( x\circ (y+\epsilon) - x\circ y  - \bigl(x^{(1)}\circ y\bigr)*\epsilon\bigr)_n\bigr|
    &\leq&
    \tfrac{1}{2} \bigl( \bigl(x^{(2)}_+ \circ (y_+ + \epsilon_+)\bigr)*\epsilon_+
    * \epsilon_+ \bigr)_n.
\end{eqnarray}
\end{prop}

\begin{pf}
For $N\geq 1$, consider new formal power series $x_N(z)=\sum_{n=1}^N
x_n z^n$, $y_N(z)=\sum_{n=1}^N y_n z^n$ and $\epsilon_N
(z)=\sum_{n=1}^{N} \epsilon_n z^n$. Since the power series $x_N\circ
(y_N+\epsilon_N)  - x_N\circ y_N$ and $x\circ (y+\epsilon)  - x\circ
y$, and $(x^{(1)}_{N,+} \circ (y_{N,+} +
\epsilon_{N,+}))*\epsilon_{N,+}$ and $(x^{(1)}_+ \circ (y_+ +
\epsilon_+))*\epsilon_+ $ have the\vspace*{1pt} same first $N$ elements, and since
$N$ is arbitrary, it is enough to prove (\ref{e:first-order}) for any
$x_N$, $y_N$ and $\epsilon_N$. Letting $f_m(u_1,\ldots,u_N)=x_m
(\sum_{k=1}^N u_k z^k  )^m$ be a real-valued function on $\Rset^N$, and
with $D f_m$ denoting the gradient of $f_m$, observe that
\begin{eqnarray}\label{expansion-aux1}
\hspace*{-30pt}\bigl(x_N\circ (y_N+\epsilon_N)  - x_N\circ y_N\bigr) (z)  &=&\sum_{m=1}^{N} x_m
 \Biggl( \Biggl(\sum_{k=1}^{N} (y_k+\epsilon_k) z^k  \Biggr)^m -
\Biggl(\sum_{k=1}^{N} y_k z^k  \Biggr)^m
\Biggr) \nonumber \\
&=& \sum_{m=1}^N \bigl( f_m(y_1+\epsilon_1,\ldots,y_N+\epsilon_N) -
f_m(y_1,\ldots,y_N) \bigr)\nonumber \\ [-8pt]\\ [-8pt]
&=& \sum_{m=1}^N  [\epsilon_1\cdots\epsilon_N]^{\mathrm{T}} D
f_m(c_{m,1},\ldots,c_{m,N})
 \nonumber \\
&=& \sum_{m=1}^N x_m m  \Biggl( \sum_{k=1}^N c_{m,k} z^k  \Biggr)^{m-1} \sum_{k=1}^N
\epsilon_k z^k,\nonumber
\end{eqnarray}
where $A^{\mathrm{T}}$ is the transpose of the matrix $A$ and
$(c_{m,1},\ldots,c_{m,N})=(1-t_m)(y_1,\ldots,y_N)+t_m(y_1+\epsilon_1,\ldots,y_N+\epsilon_N)$
for some $t_m\in[0,1]$. On the other hand, observe that
\begin{equation}\label{expansion-aux2}
\bigl(\bigl(x^{(1)}_{N,+} \circ (y_{N,+} + \epsilon_{N,+})\bigr)* \epsilon_{N,+}\bigr) (z)
= \sum_{m=1}^N x_m^+ m  \Biggl( \sum_{k=1}^N (y_k^+ + \epsilon^+_k) z^k
 \Biggr)^{m-1}  \sum_{k=1}^N \epsilon_k^+ z^k .
\end{equation}
Since $|c_{m,k}| \leq  y_k^+ + \epsilon^+_k$, it is clear that the
absolute value of the $n$th element of (\ref{expansion-aux1}) is less
than or equal to the $n$th element of (\ref{expansion-aux2}). This
completes the proof of (\ref{e:first-order}).\

The proof of (\ref{e:second-order}) is similar, but involves the
observation that
\begin{eqnarray*}
&&\bigl( x_N \circ (y_N+\epsilon_N) -
x_N\circ y_N  - \bigl(x^{(1)}_N\circ y_N\bigr)*\epsilon_N\bigr)(z)  \\
 &&\quad=\sum_{m=1}^N  \bigl( f_m(y_1+\epsilon_1,\ldots,y_N+\epsilon_N)
- f_m(y_1,\ldots,y_N)-[\epsilon_1 \cdots \epsilon_N]^{\mathrm{T}}
Df_m(y_1,\ldots,y_N)  \bigr) \\
&&\quad= \frac{1}{2}\sum_{m=1}^N [\epsilon_1 \cdots \epsilon_N]^{\mathrm{T}} D^2
f_m(d_{m,1},\ldots,d_{m,N}) [\epsilon_1 \cdots \epsilon_N] \nonumber \\
&&\quad= \frac{1}{2}\sum_{m=1}^N x_m m(m-1)  \Biggl(\sum_{k=1}^N d_{m,k} z^k
\Biggr)^{m-2}  \Biggl(\sum_{k=1}^N \epsilon_k z^k  \Biggr)^2,
\end{eqnarray*}
where $D^2 f$ denotes the Hessian of $f_m$ and
$(d_{m,1},\ldots,d_{m,N})=(1-t_m)(y_1,\ldots,y_N)+t_m(y_1+\epsilon_1,\ldots,y_N+\epsilon_N)$
with $t_m\in[0,1]$.
\end{pf}
\end{appendix}

\section*{Acknowledgements}
The authors would like to thank the anonymous referees for useful
comments and suggestions. The second author was supported in part by
NSF Grant DMS-06-08669.

\printhistory

\end{document}